# ON THE STRUCTURE OF SOLUTIONS OF ERGODIC TYPE BELLMAN EQUATION RELATED TO RISK-SENSITIVE CONTROL

By Hidehiro Kaise and Shuenn-Jyi Sheu[1]

*Nagoya University and Academia Sinica*

Bellman equations of ergodic type related to risk-sensitive control are considered. We treat the case that the nonlinear term is positive quadratic form on first-order partial derivatives of solution, which includes linear exponential quadratic Gaussian control problem. In this paper we prove that the equation in general has multiple solutions. We shall specify the set of all the classical solutions and classify the solutions by a global behavior of the diffusion process associated with the given solution. The solution associated with ergodic diffusion process plays particular role. We shall also prove the uniqueness of such solution. Furthermore, the solution which gives us ergodicity is stable under perturbation of coefficients. Finally, we have a representation result for the solution corresponding to the ergodic diffusion.

**1. Introduction.** We consider the following nonlinear partial differential equation:

$$(1.1) \quad \tfrac{1}{2}D_i(a^{ij}D_jW) + \tfrac{1}{2}\hat{a}^{ij}D_iW D_jW + b \cdot \nabla W + V = \Lambda \qquad \text{in } \mathbb{R}^N,$$

or equivalently

$$\tfrac{1}{2}a^{ij}D_{ij}W + \tfrac{1}{2}\hat{a}^{ij}D_iW D_jW + \tilde{b} \cdot \nabla W + V = \Lambda,$$

$$(1.2) \qquad\qquad\qquad \tilde{b}^i(x) \equiv b^i(x) + \tfrac{1}{2}D_j a^{ij}(x),$$

where $a(x) = [a^{ij}(x)]$, $\hat{a}(x) = [\hat{a}^{ij}(x)]$ are symmetric matrices, $b(x) = (b^1(x), \dots, b^N(x))$ is a mapping of $\mathbb{R}^N$ into $\mathbb{R}^N$, and $V(x)$ is a function on $\mathbb{R}^N$. Here we utilize the notation $D_{ij} = \partial^2/\partial x_i \partial x_j$, $D_i = \partial/\partial x_i$ and the summation convention for multiple indexes. We think of a pair $(W, \Lambda)$ of function $W(x)$

Received July 2003; revised August 2004.

[1]Supported by NSC Grant 92-2115-M-001-035.

and constant $\Lambda$ as a solution of (1.1). Equation (1.1) is called an ergodic type Bellman equation. Such equations have been treated in ergodic control problems (cf. [1]). In ergodic control problems, $\hat{a}$ is *negative-definite* and more general forms of (1.1) have been studied under rather general conditions (cf. [2]). On the other hand, (1.1) also appears in risk-sensitive problems in infinite time horizon and has been studied under certain conditions (cf. [10, 13, 14, 20]). One of the main features of (1.1) in risk-sensitive control is that $\hat{a}$ might be *indefinite*. Indeed, the following equation is studied in a risk-sensitive control problem (see [10]):

$$(1.3) \quad \frac{1}{2}a^{ij}D_{ij}W + \frac{\theta}{2}a^{ij}D_iWD_jW + \inf_{z \in Z}\{f(x,z) \cdot \nabla W + L(x,z)\} = \Lambda,$$

where $f : \mathbb{R}^N \times Z \to \mathbb{R}^N$, $L : \mathbb{R}^N \times Z \to \mathbb{R}$, $Z$ is a Borel subset in $\mathbb{R}^M$ and $\theta$ is a constant in $\mathbb{R} \backslash \{0\}$ which is called a risk-sensitive parameter. Equation (1.3) is considered to characterize a logarithm-exponential type criterion per unit time on infinite time:

$$(1.4) \qquad \qquad \Lambda = \inf_{z.} \liminf_{T \to \infty} \frac{1}{T\theta} \log E[e^{\theta \int_0^T L(X_t, z_t)\,dt}],$$

where $\{X_t\}$ is a controlled diffusion process satisfying

$$dX_t = f(X_t, z_t)\,dt + \sigma(X_t)\,dB_t, \qquad X_0 = x \in \mathbb{R}^N, \qquad a^{ij}(x) = (\sigma\sigma^T)^{ij}(x),$$

$\{B_t\}$ is standard Brownian motion and $\{z_t\}$ is a $Z$-valued process which is considered as a control. The infimum in (1.4) is taken over some class of $\{z_t\}$. In particular, if we take $f(x,z) = b(x) + C(x)z$, $L(x,z) = V(x) + (1/2)z^T S(x)z$, $Z = \mathbb{R}^M$, where $C(x)$, $S(x)$ are matrices with suitable dimension and $S(x)$ is positive-definite, then (1.3) reads

$$\frac{1}{2}a^{ij}D_{ij}W + \frac{1}{2}(\theta a - CS^{-1}C^T)^{ij}D_iWD_jW + b \cdot \nabla W + V = \Lambda.$$

Note that the sign of nonlinear term $\hat{a} = \theta a - CS^{-1}C^T$ depends on $\theta$. We also remark that the infimum in (1.3) is attained at $z = -S(x)^{-1}C(x)^T\nabla W(x)$. We are concerned with the case that $\theta a - CS^{-1}C^T$ is positive-definite since in this paper we shall study the solutions of (1.1) in the case that $\hat{a}$ is *positive-definite*. Recently, it has also become known that this case happens in some investment problems in mathematical finance (cf. [3, 8, 11, 12, 21]). However, we remark that, unlike these papers, the verification theorem will not be considered in this paper. The verification theorem is to show $\Lambda$ in (1.4) is equal to $\Lambda^*$ in Theorem 2.6 and

$$z_t^* = -S(X_t)^{-1}C(X_t)^T\nabla W^*(X_t)$$

is a feedback optimal control, $W^*$ is a solution corresponding to $\Lambda^*$ [$W^*$ is usually unique if $W^*(0) = 0$]. See also [15] for some examples from investment problems. The relation between the drift term $\hat{a}\nabla W^*$ in (1.8) for $W = W^*$ and $z_t^*$ as well as its role in the risk-sensitive control problem can be seen from the arguments in [11, 12]. The main merit of our study is to show that multiple solutions exist in general for such equations. We also provide particular solution(s) that is(are) responsible for the verification theorem. We observe that the case when $\theta a - CS^{-1}C^T$ is negative-definite can also be treated by considering the equation for $(-W)$. Therefore, according to Theorem 2.6, we have the following interesting observation. Assume $c_1 \leq a(x) \leq c_2$ and $c_1 \leq C(x)S(x)^{-1}C(x)^T \leq c_2$ for some constants $c_1, c_2 > 0$. Then for small $\theta > 0$, there is $\Lambda^*$ (depending on $\theta$) such that the above equation has solution if and only if $\Lambda \leq \Lambda^*$. For large $\theta > 0$, there is $\Lambda^*$ (depending on $\theta$) such that the above equation has solution if and only if $\Lambda \geq \Lambda^*$. In a risk-sensitive control problem, it is more interesting to assume $V(x) \to \infty$ as $|x| \to \infty$. In this case, it may happen that $\Lambda^* = \infty$ for large $\theta$. See some discussion in [20].

As we mentioned in the above, the studies of solutions for Bellman equations from an analytical point of view are considered to be fundamental to determine an optimal control. Note that solutions of (1.1) have ambiguity of additive constant, that is, if $(W, \Lambda)$ is a solution of (1.1), $W(x) + c$ still satisfies (1.1) for each constant $c$. As some examples show, it is known that (1.1) has multiple solutions even if we identify the solutions up to additive constants. So, it is important to study how we pick up a particular solution of (1.1) which gives an optimal control for the problems at hand. A common way to obtain a particular solution for ergodic type Bellman equations is to study the discounted type equations. The discounted type Bellman equation corresponding to (1.1) is as follows:

$$\tfrac{1}{2}D_i(a^{ij}D_jW_\alpha) + \tfrac{1}{2}\hat{a}^{ij}D_iW_\alpha D_jW_\alpha + b \cdot \nabla W_\alpha + V = \alpha W_\alpha.$$

$\alpha > 0$ is called a discount factor. Under certain conditions, it is shown that $W_\alpha(x) - W_\alpha(x_0)$ normalized at some point $x_0 \in \mathbb{R}^N$ and $\alpha W_\alpha$ converge to some function $W(x)$ and some constant $\Lambda$, respectively. Furthermore $(W, \Lambda)$ satisfies (1.1) (cf. [10, 13, 14]). Under the conditions including the linear exponential quadratic Gaussian (LEQG) control problem, we need to consider the case that $b(x)$ [resp. $V(x)$] is at most linearly growing (resp. quadratically growing). Under such settings, $W$ is characterized to meet some growth condition and $(W, \Lambda)$ obtained by this procedure is considered to be the right solution (cf. [13, 14]).

In the present paper we directly tackle (1.1) without the procedure using the discounted type equation under the conditions including the LEQG case. We shall specify the set of $\Lambda$ for which (1.1) has a smooth solution. Furthermore we shall characterize the set of $\Lambda$ by noting the global behavior of diffusion process which is related to some control problem.

To explain how we relate (1.1) to a control problem, we shall give a control interpretation to (1.1). Let $(\Omega, \mathcal{F}, P, \{\mathcal{F}_t\})$ be a probability space with filtration. Consider the following controlled stochastic differential equation (SDE):

$$dX_t = (\tilde{b}(X_t) + u_t)\, dt + \sigma(X_t)\, dB_t, \qquad X_0 = x \in \mathbb{R}^N, \sigma(x) \equiv a(x)^{1/2},$$

where $\{B_t\}$ is $N$-dimensional $\{\mathcal{F}_t\}$-Brownian motion and $\{u_t\}$ is an $\{\mathcal{F}_t\}$-progressively measurable process taking its value in $\mathbb{R}^N$. $\{u_t\}$ is considered as control process. We define the value function as follows:

$$v(t,x) = \sup_{u.} E_x\left[ \int_0^{T-t} (V(X_s) - \tfrac{1}{2}\hat{a}_{ij}^{-1}(X_s)u_s^i u_s^j)\, ds \right],$$

where $\hat{a}_{ij}^{-1}$ is the $(i,j)$-component in inverse of $\hat{a}$. By using the Bellman principle, we see that $v(t,x)$ satisfies the following equation formally:

$$(1.5) \quad \frac{\partial v}{\partial t} + \frac{1}{2}a^{ij}D_{ij}v + \sup_{u \in \mathbb{R}^N}\left\{ (\tilde{b}(x) + u)\cdot\nabla_x v - \frac{1}{2}\hat{a}_{ij}^{-1}u^i u^j \right\} + V = 0$$
$$\text{in } (0,T)\times\mathbb{R}^N,$$

$$(1.6) \quad v(T,x) = 0, \qquad x \in \mathbb{R}^N.$$

Since $\sup_{u \in \mathbb{R}^N}\{(\tilde{b}+u)\cdot\nabla_x v - (1/2)\hat{a}_{ij}^{-1}u^i u^j\} = (1/2)\hat{a}^{ij}D_i v D_j v + \tilde{b}\cdot\nabla_x v$, (1.5) reduces to the following:

$$\frac{\partial v}{\partial t} + \frac{1}{2}a^{ij}D_{ij}v + \frac{1}{2}\hat{a}^{ij}D_i v D_j v + \tilde{b}\cdot\nabla_x v + V = 0.$$

Note that the supremum is attained at $\bar{u}(t,x) = \hat{a}(x)\nabla_x v(t,x)$. If $-(\partial v/\partial t)(0,x)$ converges to some constant $\Lambda$ and $v(0,x) - v(0,x_0)$ normalized at some point $x_0 \in \mathbb{R}^N$ converges to some function $W(x)$ as $T \to \infty$, we have formally the following equation which we shall discuss in this paper:

$$\tfrac{1}{2}a^{ij}D_{ij}W + \tfrac{1}{2}\hat{a}^{ij}D_i W D_j W + \tilde{b}\cdot\nabla W + V = \Lambda.$$

This is considered to characterize the long-time average cost defined as following:

$$(1.7) \quad \Lambda = \sup_{u.}\limsup_{T\to\infty}\frac{1}{T}E_x\left[ \int_0^T \left( V(X_s) - \frac{1}{2}\hat{a}_{ij}^{-1}(X_s)u_s^i u_s^j \right) ds \right].$$

Following the Bellman principle, we can expect that $\bar{u}_t = \hat{a}(X_t)\nabla W(X_t)$ should be a candidate of optimal control for (1.7), where $\{X_t\}$ is defined by the controlled SDE with $u_t = \bar{u}_t = \hat{a}(X_t)\nabla W(X_t)$:

$$(1.8) \quad dX_t = (\tilde{b}(X_t) + \hat{a}\nabla W(X_t))\, dt + \sigma(X_t)\, dB_t, \qquad X_0 = x.$$

We shall study the structure of solutions of (1.1) by relating to (1.8) under conditions which include the LEQG case, that is, $b(x)$ [resp. $V(x)$] has at most linear growth (resp. quadratic growth).

The paper is organized as follows.

In Section 2 we shall specify the set of $\Lambda$ for which (1.1) has a solution under rather general conditions on $b(x)$ and $V(x)$. Indeed, it is proved that the set of $\Lambda$ is equal to closed half-line $[\Lambda^*, \infty)$ for some $\Lambda^* \in (-\infty, \infty)$.

In Section 3 we shall classify $\Lambda$ according to the global property of the diffusion process defined by (1.8). We shall prove that for $\Lambda > \Lambda^*$, the diffusion process $\{X_t\}$ in (1.8) corresponding to solution $(W, \Lambda)$ is transient and for $\Lambda = \Lambda^*$, $\{X_t\}$ is ergodic. Moreover, we shall show that solution $W(x)$ corresponding to $\Lambda^*$ is unique up to additive constant.

We show the structure of solutions in Sections 2 and 3. In Section 4 we shall consider the problem that the structures specified in Sections 2 and 3 are preserved under the perturbation on coefficients in (1.1). More precisely, consider (1.1) with $a = a_n$, $\hat{a} = \hat{a}_n$, $b = b_n$, $V = V_n$:

$$\tfrac{1}{2} D_i(a_n^{ij} D_j W_n) + \tfrac{1}{2}\hat{a}^{ij} D_i W_n D_j W_n + b \cdot \nabla W_n + V_n = \Lambda_n.$$

In similar ways to Sections 2 and 3 we can find $[\Lambda_n^*, \infty)$ for (1.1) parameterized by $n$ and solution $W_n$ corresponding to $\Lambda_n^*$ is unique. In Section 4 we mainly study the case that $a_n = a, \hat{a}_n = \hat{a}, b_n = b$, independent of $n$, and shall show that if $V_n$ converges to $V$, $\Lambda_n^*$ converges to $\Lambda^*$ and unique solution $W_n$ corresponding to $\Lambda_n^*$ converges to unique solution corresponding to $\Lambda^*$.

In Section 5 we shall study the representation for $\Lambda^*$. To obtain the representation result, we consider perturbation on $V$ and notice the dependence on $V$ for $\Lambda^*$. By using the representation, we can prove the moment condition for invariant measure of the ergodic diffusion process in (1.8) corresponding to $\Lambda^*$.

Last, we mention the connection to positive solutions of linear equations. Suppose $a^{ij}(x) = \hat{a}^{ij}(x)$, $i, j = 1, \ldots, N$. If we take the transformation $\phi(x) \equiv e^{W(x)}$ in (1.1), we have

$$(1.9) \qquad \tfrac{1}{2} a^{ij} D_{ij}\phi + \tilde{b} \cdot \nabla \phi + V\phi = \Lambda\phi.$$

Thus, in the case that $a^{ij}(x) = \hat{a}^{ij}(x)$, the study of solutions for (1.1) reduces to that of positive solutions for (1.9). We note that the structure of $\Lambda$ specified in this paper is considered to be a generalization in the theory of positive harmonic function for linear differential operators (cf. [22]). Some applications of our results to the evaluation of large time asymptotics of expectations of diffusion processes will be given in [17].

**2. Set of $\Lambda$ with solutions.** In the present section we shall consider the set of $\Lambda$ for which (1.1) has a classical solution $W$ under rather general conditions. In the next section we shall classify $\Lambda$ by following the global behavior of the diffusion process related to the solution $W$ corresponding to $\Lambda$.

We define the following set:

(2.1)   $\mathcal{A} \equiv \{\Lambda \colon \text{ there exists smooth function } W \text{ satisfying (1.1) for } \Lambda\}$.

Under the assumptions given below, we can prove that $\mathcal{A}$ has the following form for some $\Lambda^* \in (-\infty, \infty)$:

$$\mathcal{A} = [\Lambda^*, \infty).$$

For simplicity, we always assume $a^{ij}$, $\hat{a}^{ij}$, $b$, $V$ are sufficiently smooth. We shall give the following assumptions:

(A1) There exist $0 < \nu_1 < \nu_2$ such that

$$\nu_1 |\xi|^2 \le a^{ij}(x)\xi_i\xi_j \le \nu_2 |\xi|^2 \qquad \forall\, x, \xi \in \mathbb{R}^N.$$

(A2) There exist $0 < \mu_1 < \mu_2$ such that

$$\mu_1 |\xi|^2 \le \hat{a}^{ij}(x)\xi_i\xi_j \le \mu_2 |\xi|^2 \qquad \forall\, x, \xi \in \mathbb{R}^N.$$

(A3) There exists a smooth function $W_0(x)$ such that

$$\tfrac{1}{2}D_i(a^{ij}D_jW_0) + \tfrac{1}{2}\hat{a}^{ij}D_iW_0D_jW_0 + b\cdot\nabla W_0 + V \to -\infty \qquad \text{as } |x| \to \infty.$$

REMARK 2.1.   Note that it follows from (A1), (A2) that there exist $c$, $\bar{c} > 0$ such that

(2.2)          $ca(x) \le \hat{a}(x) \le \bar{c}a(x), \qquad x \in \mathbb{R}^N.$

REMARK 2.2.   In the following, we give some interesting examples for (A3).

   (a) Assume $a_{ij}(x), \hat{a}_{ij}(x)$ are bounded together with their derivatives. Assume also that there are $c_0, r_0 > 0$ such that

$$b(x) \cdot x \le -c_0|x|^2, \qquad |x| \ge r_0.$$

   (b) Assume $a_{ij}(x), \hat{a}_{ij}(x)$ are bounded together with their derivatives. Assume also that there are $c_0, r_0 > 0$ such that

$$b(x) \cdot x \ge c_0|x|^2, \qquad |x| \ge r_0.$$

   (c) Assume $V(x) \to -\infty$ as $|x| \to \infty$.

For (a), we take $W_0(x) = c|x|^2$ for small $c > 0$. For (b), we take $W_0(x) = -c|x|^2$ for small $c > 0$. For (c), we take $W_0 = 0$.

REMARK 2.3.   For the purpose of discussion in the present section, we can replace (A3) with the existence of a super solution of (1.1) for some $\Lambda$ to ensure that $\mathcal{A} \ne \varnothing$. We need (A3) to classify $\Lambda$ in the next section. In the following, we show for any $\Lambda$ and $\tilde{R} > 0$, we can construct a subsolution of (1.2) in $B_{\tilde{R}}$ with boundary value $W_0$ on $\partial B_{\tilde{R}}$, where $B_{\tilde{R}}$ is open ball with radius $\tilde{R}$ centered at 0. Indeed, we consider linear partial differential equation with Dirichlet boundary condition:

$$\tfrac{1}{2}D_i(a^{ij}D_j\tilde{W}_0) + b \cdot \nabla \tilde{W}_0 + V = \Lambda \qquad \text{in } B_{\tilde{R}},$$

$$\tilde{W}_0(x) = W_0(x) \qquad \text{on } \partial B_{\tilde{R}}.$$

Under (A1) and smoothness of coefficients, we have unique solution $\tilde{W}_0 \in C^{2,\alpha}(B_{\tilde{R}}) \cap C(\bar{B}_{\tilde{R}})$. By (A2), we have

$$\tfrac{1}{2}D_i(a^{ij}D_j\tilde{W}_0) + \tfrac{1}{2}\hat{a}^{ij}D_i\tilde{W}_0 D_j\tilde{W}_0 + b \cdot \nabla \tilde{W}_0 + V \geq \Lambda \qquad \text{in } B_{\tilde{R}}.$$

In order to see $\mathcal{A} \neq \varnothing$, consider the following Dirichlet problem:

$$(2.3) \quad \tfrac{1}{2}D_i(a^{ij}D_jW_R) + \tfrac{1}{2}\hat{a}^{ij}D_iW_R D_jW_R + b \cdot \nabla W_R + V = \Lambda \qquad \text{in } B_R,$$

$$(2.4) \qquad\qquad\qquad\qquad\qquad\qquad\qquad W_R = W_0 \qquad \text{on } \partial B_R,$$

where $B_R$ is open ball with radius $R$ centered at 0 and $W_0$ is taken from (A3). Note that (2.3) is equivalent to

$$(2.5) \quad \tfrac{1}{2}a^{ij}D_{ij}W_R + \tfrac{1}{2}\hat{a}^{ij}D_iW_R D_jW_R + \tilde{b} \cdot \nabla W_R + V = \Lambda \qquad \text{in } B_R.$$

By (A3), $W_0$ satisfies the following inequality for some $\Lambda$:

$$\tfrac{1}{2}D_i(a^{ij}D_jW_0) + \tfrac{1}{2}\hat{a}^{ij}D_iW_0 D_jW_0 + b \cdot \nabla W_0 + V \leq \Lambda \qquad \text{in } \mathbb{R}^N.$$

Also, from Remark 2.3, we see that for $\tilde{R} > R$, there exists a smooth function $\tilde{W}_0(x)$ such that

$$\tfrac{1}{2}D_i(a^{ij}D_j\tilde{W}_0) + \tfrac{1}{2}\hat{a}^{ij}D_i\tilde{W}_0 D_j\tilde{W}_0 + b \cdot \nabla \tilde{W}_0 + V \geq \Lambda \qquad \text{in } B_{\tilde{R}}.$$

Then, under (A1)–(A3), there exists $W_R \in C^{2,\alpha}(\bar{B}_R)$ satisfying (2.3) and (2.4) (cf. [18], Chapter 4, Theorem 8.4).

We need a uniform bound for $\nabla W_R$ on compact sets to obtain a solution $W$ of (1.1) by sending the radius $R$ to $\infty$. The following gradient estimate is also useful in the later discussions.

LEMMA 2.4. *Let $W_R$ be a smooth function satisfying (2.3) in $B_R$. Under* (A1) *and* (A2), *we have for each $r > 0$ and $R > 2r$*

$$(2.6) \qquad\qquad \sup_{B_r} |\nabla W_R|^2 \leq C_r + C\Lambda,$$

*where $C$ is a nonnegative constant independent of $r$ and $R$, and $C_r$ is a constant depending only on $r$.*

PROOF.  Equation (1.1) has the nonlinear term similar to those treated in [13, 14] and we can follow the same arguments to obtain the gradient estimate. However, we shall give a proof to specify the dependence of $\Lambda$.

We set $W = W_R$ for simplicity. By differentiating each side of (2.5) on $x_k$, we have

$$\begin{aligned}
(2.7) \quad & \tfrac{1}{2} D_k a^{ij} D_{ij} W + \tfrac{1}{2} a^{ij} D_{ijk} W + \tfrac{1}{2} D_k \hat{a}^{ij} D_i W D_j W \\
& + \hat{a}^{ij} D_i W D_{jk} W + D_k \tilde{b}^i D_i W + \tilde{b}^i D_{ik} W + D_k V = 0.
\end{aligned}$$

Let us set $G \equiv (1/2) \sum_k (D_k W)^2$. Then, using (2.7)

$$\begin{aligned}
(2.8) \quad & -\tfrac{1}{2} a^{ij} D_{ij} G - \hat{a}^{ij} D_i W D_j G - \tilde{b}^i D_i G \\
& = -\tfrac{1}{2} a^{ij} D_k W D_{ijk} W - \tfrac{1}{2} a^{ij} D_{ki} W D_{kj} W \\
& \quad - \hat{a}^{ij} D_i W D_k W D_{jk} W - \tilde{b}^i D_k W D_{ik} W \\
& = \tfrac{1}{2} D_k a^{ij} D_k W D_{ij} W + \tfrac{1}{2} D_k \hat{a}^{ij} D_i W D_j W D_k W \\
& \quad + D_k \tilde{b}^i D_i W D_k W + D_k V D_k W - \tfrac{1}{2} a^{ij} D_{ki} W D_{kj} W.
\end{aligned}$$

We note the second-order derivative terms. Then, we have

$$\begin{aligned}
\text{RHS of (2.8)} & \leq \frac{1}{4\delta} \left( \sum_{i,j} |D a^{ij}|^2 \right) |DW|^2 + \frac{\delta}{4} |D^2 W|^2 \\
& \quad + \tfrac{1}{2} D_k \hat{a}^{ij} D_i W D_j W D_k W + D_k \tilde{b}^i D_i W D_k W \\
& \quad + D_k V D_k W - \tfrac{1}{4} a^{ij} D_{ki} W D_{kj} W - \tfrac{1}{4} a^{ij} D_{ki} W D_{kj} W \\
& \leq \frac{1}{4\delta} \left( \sum_{i,j} |D a^{ij}|^2 \right) |DW|^2 + \tfrac{1}{2} D_k \hat{a}^{ij} D_i W D_j W D_k W \\
& \quad + D_k \tilde{b}^i D_i W D_k W + D_k V D_k W - \tfrac{1}{4} a^{ij} D_{ki} W D_{kj} W,
\end{aligned}$$

where $\delta > 0$ is a small constant. Indeed, we can take $\delta$ satisfying $\delta < \nu_1$. From matrix inequality $(\operatorname{tr} AB)^2 \leq N \nu_2 (\operatorname{tr} AB^2)$ where $A, B$ are $N \times N$-symmetric matrices, $A$ is nonnegative-definite and $\nu_2$ is the maximum eigenvalue of $A$, we finally obtain the following inequality:

$$\begin{aligned}
(2.9) \quad & -\frac{1}{2} a^{ij} D_{ij} G - \hat{a}^{ij} D_i W D_j G - \tilde{b}^i D_i G \\
& \leq C_r |DW| + C_r |DW|^2 + C_r |DW|^3 - \frac{1}{4N\nu_2} (a^{ij} D_{ij} W)^2 \qquad \text{in } B_{2r}.
\end{aligned}$$

Here and in the proof below, we suppose that $C_r$ is constant depending only on $r$ and $C$ is nonnegative constant independent of $r$ and $R$.

Fix arbitrary $\xi \in B_r$ and take a cut-off function $\varphi \in C_0^\infty(\mathbb{R}^N)$ satisfying the following:

$$(2.10) \quad 0 \le \varphi \le 1 \text{ in } \mathbb{R}^N, \qquad \varphi(\xi) = 1, \qquad \varphi \equiv 0 \text{ in } B_r(\xi)^c,$$
$$|\nabla \varphi| \le C\varphi^{1/2}, \qquad |D^2\varphi| \le C,$$

where $B_r(\xi)$ is open ball with radius $r$ centered at $\xi$. Let $x_0$ be a maximum point of $\varphi G$ in $\bar{B}_r(\xi)$. By the maximum principle, we can see

$$\begin{aligned}
(2.11) \quad 0 &\le -\tfrac{1}{2}a^{ij}D_{ij}(\varphi G) - \hat{a}^{ij}D_iW D_j(\varphi G) - \tilde{b}^i D_i(\varphi G) \\
&= \varphi\{-\tfrac{1}{2}a^{ij}D_{ij}G - \hat{a}^{ij}D_iW D_jG - \tilde{b}^i D_iG\} \\
&\quad -\tfrac{1}{2}a^{ij}(D_{ij}\varphi)G - a^{ij}D_i\varphi D_jG - \hat{a}^{ij}D_j\varphi(D_iW)G - \tilde{b}^i(D_i\varphi)G \\
&\le \varphi\{-\tfrac{1}{2}a^{ij}D_{ij}G - \hat{a}^{ij}D_iW D_jG - \tilde{b}^i D_iG\} \\
&\quad + C_rG + C\varphi^{1/2}G^{3/2} \qquad \text{at } x_0,
\end{aligned}$$

where we used $0 = D(\varphi G) = GD\varphi + \varphi DG$ and (2.10). From (2.5) and (2.9), it is implied that

$$\begin{aligned}
(2.12) \quad \text{RHS of (2.11)} &\le \varphi\Big\{C_rG^{1/2} + C_rG + C_rG^{3/2} - \frac{1}{4N\nu_2}(a^{ij}D_{ij}W)^2\Big\} \\
&\quad + C_rG + C\varphi^{1/2}G^{3/2} \\
&= \varphi\Big\{C_rG^{1/2} + C_rG + C_rG^{3/2} \\
&\qquad - \frac{1}{N\nu_2}\Big(-\tfrac{1}{2}\hat{a}^{ij}D_iW D_jW - \tilde{b}^i D_iW - V + \Lambda\Big)^2\Big\} \\
&\quad + C_rG + C\varphi^{1/2}G^{3/2} \qquad \text{at } x_0.
\end{aligned}$$

By (A2), the following inequalities hold for some positive constant $\kappa$ which depends on $\mu_1$,

$$\begin{aligned}
(2.13) \quad &-\tfrac{1}{2}\hat{a}^{ij}D_iW D_jW - \tilde{b}^i D_iW - V + \Lambda \\
&\le -\frac{\mu_1}{2}|DW|^2 + C_r|DW| - V + \Lambda \\
&\le -\kappa|DW|^2 + C_r - V + \Lambda.
\end{aligned}$$

In the case that $-\kappa|DW|^2 - V + \Lambda \ge 0$ at $x_0$, we have

$$\kappa|DW|^2(x_0) \le C_r - V(x_0) + \Lambda \le C_r + \Lambda,$$

where we used $x_0 \in B_{2r}$. Since

$$\tfrac{1}{2}|DW|^2(\xi) = \tfrac{1}{2}|DW|^2(\xi)\varphi(\xi) \le G(x_0)\varphi(x_0),$$

we obtain the following gradient estimate at $\xi$:

$$\kappa |DW|^2(\xi) \leq C_r + \Lambda.$$

We next consider the case that

$$-\kappa |DW|^2 + C_r - V + \Lambda \leq 0 \qquad \text{at } x_0.$$

By (2.13),

RHS of (2.12)

$$\leq \varphi \left\{ C_r G^{1/2} + C_r G + C_r G^{3/2} - \frac{1}{N\nu_2}(-\kappa |DW|^2 + C_r - V + \Lambda)^2 \right\}$$

(2.14) $\qquad + C_r G + C\varphi^{1/2} G^{3/2}$

$$\leq \varphi \left\{ C_r G^{1/2} + C_r G + C_r G^{3/2} - \frac{4\kappa^2}{N\nu_2} G^2 + \frac{4\kappa}{N\nu_2} G(C_r - V + \Lambda) \right\}$$

$$+ C_r G + C\varphi^{1/2} G^{3/2}.$$

If $C_r - V + \Lambda \geq \kappa G(x_0)/4$ or $C_r \geq G(x_0)$ we have the bound $|DW|^2(\xi) \leq C_r + C\Lambda$ in the same way as the above case. We shall consider the case that $C_r - V + \Lambda \leq \kappa G(x_0)/4$ and $C_r \leq G(x_0)$. Then, from (2.14), we have

$$0 \leq \varphi \left\{ G^{3/2} + C_r G + C_r G^{3/2} - \frac{4\kappa^2}{N\nu_2} G^2 + \frac{2\kappa^2}{N\nu_2} G^2 \right\}$$

$$+ C_r G + C\varphi^{1/2} G^{3/2}$$

$$\leq -C_1 \varphi G^2 + C_2 \varphi^{1/2} G^{3/2} + C_3 C_r G$$

$$\equiv -C_1 \varphi G^2 + C_2 \varphi^{1/2} G^{3/2} + \tilde{C}_3 G \qquad \text{at } x_0, \tilde{C}_3 \equiv C_3 C_r,$$

where $C_1$, $C_2$, $C_3$ are positive constants independent of $r$, $R$ and $\Lambda$. By setting $X \equiv \varphi^{1/2}(x_0) G^{1/2}(x_0)$, we have

$$0 \leq -C_1 X^2 + C_2 X + \tilde{C}_3.$$

Therefore, we have

$$X^2 = \varphi G(x_0) \leq \frac{C_2^2}{C_1^2} + \frac{2\tilde{C}_3}{C_1} \leq \frac{C_2^2}{C_1^2} + \frac{2C_3 C_r}{C_1}.$$

Since $(1/2)|DW|^2(\xi) = (1/2)|DW|^2(\xi)\varphi^2(\xi) \leq G(x_0)\varphi(x_0)$, we obtain the bound for $|DW|(\xi)$.  $\square$

REMARK 2.5. Under some growth conditions for coefficients of (1.2), we can obtain growth order for gradient of solutions. For instance, besides (A1), (A2), suppose the following growth conditions: $Da^{ij}(x), D\hat{a}^{ij}(x)$ are bounded and there exist $c_1, c_2 > 0$ and $m \geq 1$ such that

$$|b(x)| \leq c_1(1 + |x|^m), \qquad |Db(x)| \leq c_1(1 + |x|^{m-1}),$$
$$|V(x)| \leq c_2(1 + |x|^{2m}), \qquad |DV(x)| \leq c_2(1 + |x|^{2m-1}).$$

Then, in Lemma 2.4, we can take $C_r = C(1 + r^{2m})$, that is,

$$\sup_{B_r} |\nabla W_R|^2 \leq C(1 + r^{2m} + \Lambda), \qquad 0 < 2r < R.$$

We may normalize $W_R$ as $W_R(0) = 0$ because (1.1) does not include a zeroth term on $W_R$. Then, from Lemma 2.4, there exists $W \in C(\mathbb{R}^N)$ such that $W_R$ converges to $W$ on each compact set as $R \to \infty$ by taking a subsequence if necessary. Also, since $\{W_R\}_{R > 2r}$ is bounded in $H^1(B_r)$ by Lemma 2.4, $W_R$ converges to $W$ $L^2_{\text{loc}}$-strongly and $H^1_{\text{loc}}$-weakly. Furthermore, we can see that $\nabla W_R$ converges $L^2_{\text{loc}}$-strongly in a similar way to Lemma 2.8 in [14] and Section 1.4 in [20].

We rewrite (2.3), (2.4) in integral form:

$$-\tfrac{1}{2} \int a^{ij} D_i W_R D_j \varphi \, dx + \tfrac{1}{2} \int \hat{a}^{ij} D_i W_R D_j W_R \varphi \, dx$$
$$+ \int b \cdot \nabla W_R \varphi \, dx + \int V \varphi \, dx = \int \Lambda \varphi \, dx, \qquad \varphi \in C_0^\infty(B_R).$$

Fix $r > 0$. Since $W_R$ converges to $W$ $H^1_{\text{loc}}$-strongly, we obtain the following by sending $R$ to $\infty$:

$$-\tfrac{1}{2} \int a^{ij} D_i W D_j \varphi \, dx + \tfrac{1}{2} \int \hat{a}^{ij} D_i W D_j W \varphi \, dx$$
$$+ \int b \cdot \nabla W \varphi \, dx + \int V \varphi \, dx = \int \Lambda \varphi \, dx, \qquad \varphi \in C_0^\infty(B_r), r > 0.$$

Owing to the regularity theorem of elliptic equations and the imbedding theorem, we have $W$ as a classical solution of (1.1). Therefore, we have proved that $\mathcal{A} \neq \varnothing$.

We shall state and prove the form of the set of $\Lambda$.

THEOREM 2.6.   *Under the assumptions* (A1)–(A3), *there exists* $\Lambda^* \in (-\infty, \infty)$ *such that* $\mathcal{A} = [\Lambda^*, \infty)$.

PROOF.   In order to show $\inf \mathcal{A} > -\infty$, we suppose $\inf \mathcal{A} = -\infty$, that is, there exists $\{\Lambda_n\} \subset \mathcal{A}$ such that $\Lambda_n$ tends to $-\infty$ as $n \to \infty$. Let $W_n$ be a solution of (1.1) corresponding to $\Lambda_n$. Then, by the integral form of (1.1), we have

$$
\begin{aligned}
(2.15) \quad &-\tfrac{1}{2}\int a^{ij}D_i W_n D_j\varphi\,dx + \tfrac{1}{2}\int \hat{a}^{ij}D_i W_n D_j W_n\,\varphi\,dx \\
&\qquad + \int b\cdot\nabla W_n\,\varphi\,dx + \int V\varphi\,dx = \int \Lambda_n\varphi\,dx, \qquad \varphi\in C_0^\infty(\mathbb{R}^N).
\end{aligned}
$$

Take $\varphi\in C_0^\infty(\mathbb{R}^N)$ such that $\int\varphi\,dx\neq 0$. Since $\{\Lambda_n\}$ is bounded from above, it is implied from Lemma 2.4 that

$$
(2.16) \qquad \sup_{B_r}|\nabla W_n|\le C_r,
$$

where $C_r$ is a constant independent of $n$ and $r$ is taken such that $\sup\varphi\subset B_r$. Therefore, the left-hand side of (2.15) is bounded on $n$. On the other hand, the right-hand side of (2.15) is unbounded because of the assumption which we made above. This leads to a contradiction.

We shall next prove if $\tilde\Lambda\in\mathcal{A}$, then $[\tilde\Lambda,\infty)\subset\mathcal{A}$. Let $\tilde W$ be a solution corresponding to $\tilde\Lambda$. For arbitrary $\Lambda\ge\tilde\Lambda$, we have

$$
(2.17) \quad \tfrac{1}{2}D_i(a^{ij}D_j\tilde W) + \tfrac{1}{2}\hat{a}^{ij}D_i\tilde W D_j\tilde W + b\cdot\nabla\tilde W + V = \tilde\Lambda\le\Lambda \qquad \text{in } \mathbb{R}^N.
$$

By Remark 2.3, for $\tilde R>R$, there exists $\tilde W_0$ such that

$$
(2.18) \quad \tfrac{1}{2}D_i(a^{ij}D_j\tilde W_0) + \tfrac{1}{2}\hat{a}^{ij}D_i\tilde W_0 D_j\tilde W_0 + b\cdot\nabla\tilde W_0 + V \ge\Lambda \qquad \text{in } B_{\tilde R}.
$$

Consider the Dirichlet problem (2.3) with boundary condition $W_R=\tilde W_0$ on $\partial B_R$. From (2.17), (2.18), the existence of a classical solution for this Dirichlet problem is guaranteed by Theorem 8.4, Chapter 4 in [18]. In the same manner as that right after the proof of Lemma 2.4, we can see that there exists a smooth function $W$ satisfying (1.1) for $\Lambda$.

We shall prove that $\Lambda^*\equiv\inf\mathcal{A}$ actually belongs to $\mathcal{A}$. $\{\Lambda_n\}$ is a sequence in $\mathcal{A}$ such that $\Lambda_n\to\Lambda^*$ and $W_n$ is a solution of (1.1) corresponding to $\Lambda_n$ normalized as $W_n(0)=0$. Then, $W_n$ satisfies (2.15). Since $\{\Lambda_n\}$ is bounded, it follows from Lemma 2.4 that (2.16) holds for some constant $C_r$ independent of $n$. Following the same way as the discussion after Lemma 2.4, we can see that a sequence of $W_n$ converges to $W^*\in C(\mathbb{R}^N)$ uniformly on compact sets and $H^1_{\mathrm{loc}}$-strongly. By taking a limit in (2.15) as $n\to\infty$, we have

$$
\begin{aligned}
&-\tfrac{1}{2}\int a^{ij}D_i W^* D_j\varphi\,dx + \tfrac{1}{2}\int \hat{a}D_i W^* D_j W^*\,\varphi\,dx \\
&\qquad + \int b\cdot\nabla W^*\varphi\,dx + \int V\varphi\,dx = \int \Lambda^*\varphi\,dx \qquad \forall\varphi\in C_0^\infty(\mathbb{R}^N).
\end{aligned}
$$

Therefore, the existence of a classical solution $W^*$ of (1.1) for $\Lambda^*$ follows from the regularity theorem of elliptic equations (see Theorems 5.1, 6.3, Chapter 4 in [18]). $\quad\square$

### 3. Classification of solutions.

3.1. *Transience and ergodicity of diffusion processes.* In the last section we proved that the set of $\Lambda$ for which (1.1) has a smooth solution is $\mathcal{A} = [\Lambda^*, \infty)$ for some $\Lambda^* \in (-\infty, \infty)$. In the present section we shall study the classification of $\Lambda$ by global behavior of $\{X_t\}$ defined by (1.8).

Let $(\Omega, \mathcal{F}, P, \{\mathcal{F}_t\})$ be a filtered probability space on which $N$-dimensional Brownian motion $\{B_t\}$ is defined. For given $\Lambda \in [\Lambda^*, \infty)$, consider the SDE:

$$(3.1) \qquad dX_t = (\tilde{b}(X_t) + \hat{a}\nabla W(X_t))\,dt + \sigma(X_t)\,dB_t, \qquad X_0 = x,$$

where $W(x)$ is a solution of (1.1) corresponding to $\Lambda$. We shall classify $\Lambda$ according to the global properties of $\{X_t\}$. More precisely, we shall prove that for $\Lambda > \Lambda^*$, $\{X_t\}$ is transient and for $\Lambda = \Lambda^*$, $\{X_t\}$ is ergodic. Note that solution of (3.1) might explode in finite time.

We shall next discuss transience of $\{X_t\}$ for $\Lambda \in (\Lambda^*, \infty)$. We introduce the operator associated to solution $(W, \Lambda)$ of (1.1):

$$T_t^{W,\Lambda} f(x) \equiv E_x[f(X_t); t < \tau_\infty], \qquad f \in C_0(\mathbb{R}^N),$$

$$\tau_n \equiv \inf\{t; X_t \notin B_n(0)\}, \qquad \tau_\infty \equiv \lim_{n \to \infty} \tau_n,$$

where $\{X_t\}$ is a solution of (3.1) up to $t < \tau_\infty$ corresponding to $(W, \Lambda)$.

LEMMA 3.1. *Under* (A1)–(A3), *the following inequality holds for each solution $(W, \Lambda)$ of* (1.1):

$$T_t^{W,\Lambda} f(x) \leq k(x)e^{-c(\Lambda - \Lambda^*)t}, \qquad f \in C_0(\mathbb{R}^N), f \geq 0,$$

*where $c$ is in Remark* 2.1 *and $k(x)$ is a constant depending only on $x$.*

PROOF. Let $W^*$ be a solution of (1.1) corresponding to $\Lambda^*$. We set $W_c \equiv cW$, $W_c^* \equiv cW^*$, where $c > 0$ is taken from Remark 2.1. Then, we have from (1.2)

$$(3.2) \qquad \frac{1}{2}a^{ij}D_{ij}W_c + \frac{1}{2c}\hat{a}^{ij}D_i W_c D_j W_c + \tilde{b} \cdot \nabla W_c + cV = c\Lambda,$$

$$(3.3) \qquad \frac{1}{2}a^{ij}D_{ij}W_c^* + \frac{1}{2c}\hat{a}^{ij}D_i W_c^* D_j W_c^* + \tilde{b} \cdot \nabla W_c^* + cV = c\Lambda^*.$$

Subtracting (3.3) from (3.2),

$$\frac{1}{2}a^{ij}D_{ij}(W_c - W_c^*) + (\tilde{b} + \hat{a}\nabla W^*) \cdot \nabla(W_c - W_c^*)$$

$$+ \frac{1}{2c}\hat{a}\nabla(W_c - W_c^*) \cdot \nabla(W_c - W_c^*) = c(\Lambda - \Lambda^*).$$

Setting $\bar{W} \equiv W_c - W_c^*$, we have

$$(3.4) \quad \frac{1}{2}a^{ij}D_{ij}\bar{W} + (\tilde{b} + \hat{a}\nabla W^*) \cdot \nabla\bar{W} + \frac{1}{2c}\hat{a}\nabla\bar{W} \cdot \nabla\bar{W} = c(\Lambda - \Lambda^*).$$

We consider (3.1) and rewrite this as follows:

$$\begin{aligned}
dX_t &= (\tilde{b}(X_t) + \hat{a}\nabla W(X_t))\,dt + \sigma(X_t)\,dB_t \\
&= (\tilde{b}(X_t) + \hat{a}\nabla W(X_t))\,dt - a\nabla\bar{W}(X_t)\,dt + \sigma(X_t)\,dB_t + a\nabla\bar{W}(X_t)\,dt \\
&= (\tilde{b}(X_t) + \hat{a}\nabla W^*(X_t))\,dt + \left(\frac{1}{c}\hat{a}\nabla\bar{W}(X_t) - a\nabla\bar{W}(X_t)\right)dt + \sigma(X_t)\,d\tilde{B}_t,
\end{aligned}$$

where

$$(3.5) \qquad \tilde{B}_s = B_s + \int_0^s \sigma\nabla\bar{W}(X_r)\,dr, \qquad s < \tau_\infty.$$

Define measure $\tilde{P}$ on $\mathcal{F}_t^{(n)} = \mathcal{F}_{t\wedge\tau_n}$:

$$\frac{d\tilde{P}}{dP}\bigg|_{\mathcal{F}_t^{(n)}} = \exp\left[-\int_0^t \sigma\nabla\bar{W}(X_s)\mathbb{1}_{\{s\leq\tau_n\}}\,dB_s - \frac{1}{2}\int_0^t a\nabla\bar{W}\cdot\nabla\bar{W}(X_s)\mathbb{1}_{\{s\leq\tau_n\}}\,ds\right].$$

Then, $\tilde{P}$ is probability measure and $\langle\tilde{B}_{\cdot\wedge\tau_n}\rangle_s = \langle B_{\cdot\wedge\tau_n}\rangle_s = s\wedge\tau_n$. By (3.5), we have

$$\begin{aligned}
&E_x[f(X_t); t < \tau_n] \\
(3.6) \quad &= \tilde{E}_x[f(X_t)e^{\int_0^t \sigma\nabla\bar{W}(X_s)\mathbb{1}_{\{s\leq\tau_n\}}\,dB_s + (1/2)\int_0^t a\nabla\bar{W}\cdot\nabla\bar{W}(X_s)\mathbb{1}_{\{s\leq\tau_n\}}\,ds}; t < \tau_n] \\
&= \tilde{E}_x[f(X_t)e^{\int_0^t \sigma\nabla\bar{W}(X_s)\mathbb{1}_{\{s\leq\tau_n\}}\,d\tilde{B}_s - (1/2)\int_0^t a\nabla\bar{W}\cdot\nabla\bar{W}(X_s)\mathbb{1}_{\{s\leq\tau_n\}}\,ds}; t < \tau_n],
\end{aligned}$$

where $\tilde{E}_x$ denotes expectation with respect to $\tilde{P}$. Applying the Itô formula to $\bar{W}(X_t)$,

$$\begin{aligned}
d\bar{W}(X_t) &= \nabla\bar{W} \cdot \left(\tilde{b} + \hat{a}\nabla W^* + \frac{1}{c}\hat{a}\nabla\bar{W} - a\nabla\bar{W}\right)(X_t)\,dt \\
&\quad + \frac{1}{2}a^{ij}D_{ij}\bar{W}(X_t)\,dt + \sigma\nabla\bar{W}(X_t)\,d\tilde{B}_t \\
&= \left(\frac{1}{2}a^{ij}D_{ij}\bar{W} + (\tilde{b} + \hat{a}\nabla W^*) \cdot \nabla\bar{W}\right)(X_t)\,dt \\
&\quad + \left(\frac{1}{c}\hat{a}\nabla\bar{W} \cdot \nabla\bar{W} - a\nabla\bar{W} \cdot \nabla\bar{W}\right)(X_t)\,dt + \sigma\nabla\bar{W}(X_t)\,d\tilde{B}_t \\
&= \left(-\frac{1}{2c}\hat{a}\nabla\bar{W} \cdot \nabla\bar{W} + c(\Lambda - \Lambda^*)\right)(X_t)\,dt \\
&\quad + \left(\frac{1}{c}\hat{a}\nabla\bar{W} \cdot \nabla\bar{W} - a\nabla\bar{W} \cdot \nabla\bar{W}\right)(X_t)\,dt + \sigma\nabla\bar{W}(X_t)\,d\tilde{B}_t
\end{aligned}$$

(3.7)

$$= \sigma \nabla \bar{W}(X_t) \, d\tilde{B}_t - \frac{1}{2} a \nabla \bar{W} \cdot \nabla \bar{W}(X_t) \, dt$$

$$+ \frac{1}{2} \left( \frac{1}{c} \hat{a} - a \right) \nabla \bar{W} \cdot \nabla \bar{W}(X_t) \, dt + c(\Lambda - \Lambda^*) \, dt.$$

Here we used (3.4). Then, by (3.6) and (3.7), we have

$$E_x[f(X_t); t < \tau_n]$$

$$= \tilde{E}_x[f(X_t)e^{-c(\Lambda - \Lambda^*)t + \bar{W}(X_t) - \bar{W}(x) + (1/2) \int_0^t (a - (1/c)\hat{a}) \nabla \bar{W} \cdot \nabla \bar{W}(X_s) \mathbb{1}_{\{s \le \tau_n\}} \, ds};$$

$$t < \tau_n]$$

$$\le \|f\|_\infty e^{\sup\{\bar{W}(y) - \bar{W}(x); y \in \operatorname{supp} f\}}$$

$$\times e^{-c(\Lambda - \Lambda^*)t} \tilde{E}_x[e^{(1/2) \int_0^t (a - (1/c)\hat{a}) \nabla \bar{W} \cdot \nabla \bar{W}(X_s) \mathbb{1}_{\{s \le \tau_n\}} \, ds}; t < \tau_n].$$

Since $ca(x) \le \hat{a}(x)$, we have

$$E_x[f(X_t); t < \tau_n] \le k(x)e^{-c(\Lambda - \Lambda^*)t},$$

$$k(x) = \|f\|_\infty \exp\left( \sup_{y \in \operatorname{supp} f} (\bar{W}(y) - \bar{W}(x)) \right).$$

Taking the limit as $n \to \infty$, we obtain

$$E_x[f(X_t); t < \tau_\infty] \le k(x)e^{-c(\Lambda - \Lambda^*)t}. \qquad \square$$

Now we have the result on transience.

THEOREM 3.2. *Let* $(W, \Lambda)$ *be a solution of* (1.1) *and* $\{X_t\}$ *be a solution of* (3.1) *corresponding to* $(W, \Lambda)$. *If* (A1)–(A3) *hold, then for* $\Lambda > \Lambda^*$, $\{X_t\}$ *is transient.*

PROOF. Let $f \in C_0(\mathbb{R}^N)$ and $f \ge 0$. Since $\Lambda > \Lambda^*$, we can see that by Lemma 3.1,

$$\int_0^\infty T_t^{W, \Lambda} f(x) \, dt < \infty.$$

Therefore, $\{X_t\}$ is transient. $\square$

We proved that for $\Lambda > \Lambda^*$, $\{X_t\}$ defined by (3.1) is transient. We next show that if $\Lambda = \Lambda^*$, the corresponding diffusion process $\{X_t^*\}$ satisfying (3.1) is ergodic.

We have to show the following proposition.

PROPOSITION 3.3. *Let $(W, \Lambda)$ be a solution of* (1.1) *and let $\{X_t\}$ be the corresponding diffusion process defined by* (3.1). *Assume* (A1)–(A3). *If $\{X_t\}$ is transient, then there exists $\alpha > 0$ such that*

$$T_t^{W,\Lambda} f(x) \leq C(x) e^{-\alpha t}, \qquad f \in C_0(\mathbb{R}^N), f \geq 0, x \in \mathbb{R}^N,$$

*where $C(x)$ is a constant independent of $t$, but depending on $x$.*

We prepare several lemmas to prove the above proposition.

Let $(W, \Lambda)$ be a solution of (1.1) and $\{X_t\}$ be a solution of (3.1). We define occupation measure for $\{X_t\}$ on $\{t < \tau_\infty\}$ as follows:

$$\mu_t(B) \equiv \frac{1}{t} \int_0^t \mathbb{1}_B(X_s)\, ds, \qquad B \in \mathcal{B}(\mathbb{R}^N), t < \tau_\infty,$$

where $\mathcal{B}(\mathbb{R}^N)$ is the Borel $\sigma$-field on $\mathbb{R}^N$. Let $\mathcal{M}_1(\mathbb{R}^N)$ be the set of probability measures on $\mathcal{B}(\mathbb{R}^N)$. We think of $\mathcal{M}_1(\mathbb{R}^N)$ as the topological vector space with topology compatible to weak convergence. Note that $\mu_t \in \mathcal{M}_1(\mathbb{R}^N)$ on $\{t < \tau_\infty\}$.

The following lemma on large deviation type estimate is useful.

LEMMA 3.4. *Let $\{X_t\}$ be a solution of* (3.1). *Then, the following estimate holds:*

$$\limsup_{t \to \infty} \frac{1}{t} \log P[\mu_t \in \mathcal{K}, t < \tau_\infty] \leq -\inf_{\mu \in \mathcal{K}} I^W(\mu)$$

(3.8)

$$\mathcal{K} \text{ is compact set in } \mathcal{M}_1(\mathbb{R}^N).$$

$I^W(\mu)$ *is defined as follows*:

$$I^W(\mu) \equiv -\inf_{u \in \mathcal{U}} \int \frac{Lu}{u}(x)\mu(dx), \qquad L \equiv \tfrac{1}{2} a^{ij} D_{ij} + (\tilde{b} + \hat{a} \nabla W) \cdot \nabla,$$

$$\mathcal{U} \equiv \{u \in C^2(\mathbb{R}^N) : u(x) > 0 \text{ for all } x, Lu/u \text{ is bounded above}\}.$$

Note that $I^W(\mu)$ takes values on $[0, \infty]$ and is convex, lower semi-continuous on $\mathcal{M}_1(\mathbb{R}^N)$. This type of estimate is well known in large deviation theory. As noted in [5], even if the state space of $\{X_t\}$ is not compact, (3.8) holds for compact set $\mathcal{K}$ (cf. comments in Section 7, page 440 of [5] and see the proof in Section 2.2 of [4] for Brownian motion).

LEMMA 3.5. *If $I^W(\mu^*) = 0$ for some $\mu^* \in \mathcal{M}_1(\mathbb{R}^N)$, then diffusion process $\{X_t\}$ defined in* (3.1) *does not explode in finite time.*

PROOF. From assumption $I^W(\mu^*) = 0$, it is implied that

$$\int \frac{Lu}{u}\, d\mu^* \geq 0 \qquad \forall u \in \mathcal{U}.$$

(3.9)

For $u \in C_0^\infty(\mathbb{R}^N)$, $u \geq 0$, and constant $c > 0$,

$$\frac{d}{dt} \int \log \frac{T_t^{W,\Lambda} u + c}{u + c} \, d\mu^* = \int \frac{T_t^{W,\Lambda} L u}{T_t^{W,\Lambda} u + c} \, d\mu^*$$

$$= \int \frac{L T_t^{W,\Lambda} u}{T_t^{W,\Lambda} u + c} \, d\mu^* = \int \frac{L(T_t^{W,\Lambda} u + c)}{T_t^{W,\Lambda} u + c} \, d\mu^*.$$

Since $T^{W,\Lambda} u + c \in \mathcal{U}$, we have

$$\frac{d}{dt} \int \log \frac{T_t^{W,\Lambda} u + c}{u + c} \, d\mu^* \geq 0 \qquad \forall \, t.$$

Thus, we can see that

(3.10)
$$\int \log \frac{T_t^{W,\Lambda} u + c}{u + c} \, d\mu^* \geq \int \log \frac{T_0^{W,\Lambda} u + c}{u + c} \, d\mu^*$$

$$= \int \log \frac{u + c}{u + c} \, d\mu^* = 0.$$

Let $\{\phi_n(x)\}_{n=1}^\infty$ be a sequence in $C_0^\infty(\mathbb{R}^N)$ such that $0 \leq \phi_n(x) \leq 1$ and $\phi_n(x) \uparrow 1$ as $n \to \infty$. If we take $u = \phi_n$ in (3.10), then we have

$$\int \log \frac{T_t^{W,\Lambda} \phi_n + c}{\phi_n + c} \, d\mu^* \geq 0.$$

Since $T_t^{W,\Lambda} \phi_n \leq T_t^{W,\Lambda} 1$,

$$\int \log \frac{T_t^{W,\Lambda} 1 + c}{\phi_n + c} \, d\mu^* \geq 0.$$

By taking the limit as $n \to \infty$, we obtain

$$\int \log \frac{T_t^{W,\Lambda} 1 + c}{1 + c} \, d\mu^* \geq 0.$$

Noting that $(T_t^{W,\Lambda} 1 + c)/(1 + c) \leq 1$, we can see that

$$T_t^{W,\Lambda} 1 = 1, \qquad \mu^*\text{-a.s.}$$

Since the diffusion process is nondegenerate [see (A1)],

$$T_t^{W,\Lambda} 1(x) = 1 \qquad \forall \, x \in \mathbb{R}^N, \forall \, t \geq 0.$$

Finally, as $t \to \infty$, we have

$$P_x[\tau_\infty = \infty] = \lim_{t \to \infty} P_x[t < \tau_\infty] = \lim_{t \to \infty} T_t^{W,\Lambda} 1(x) = 1. \qquad \square$$

LEMMA 3.6. *Let $\{X_t\}$ be a solution of (3.1). If $I^W(\mu^*) = 0$, then $\mu^*$ is an invariant measure for $\{X_t\}$.*

Proof. Since $I^W(\mu^*) = -\inf_{u \in \mathcal{U}} \int (Lu/u)(x)\mu^*(dx) = 0$,

$$\int \frac{Lu}{u}(x)\mu^*(dx) \geq 0 \qquad \forall u \in \mathcal{U}.$$

Setting $w = \log u$, we have

$$(3.11) \qquad \int (Lw(x) + \tfrac{1}{2}a\nabla w \cdot \nabla w(x))\mu^*(dx) \geq 0, \qquad u = e^w \in \mathcal{U}.$$

Denote $\tilde{\mathcal{U}}$ by

$$\tilde{\mathcal{U}} = \{u \in C^2(\mathbb{R}^N); \exists R > r > 0 \text{ s.t. } r \leq u(x) \leq R,$$
$$Du, D^2u \text{ have compact support}\}.$$

Note that $\tilde{\mathcal{U}} \subset \mathcal{U}$. It is easy to see that if $u = e^w \in \tilde{\mathcal{U}}$, then $u_\lambda \equiv e^{\lambda w} \in \tilde{\mathcal{U}}$ for $\lambda \in \mathbb{R}$. Therefore, applying $\lambda w$ in (3.11) instead of $w$,

$$\int \left(Lw(x) + \frac{\lambda}{2}a\nabla w \cdot \nabla w(x)\right)\mu^*(dx) \geq 0, \qquad u = e^w \in \tilde{\mathcal{U}}, \lambda > 0.$$

Taking the limit as $\lambda \to 0$, we have

$$\int Lw(x)\mu^*(dx) \geq 0, \qquad u = e^w \in \tilde{\mathcal{U}}.$$

Since $u = e^w \in \tilde{\mathcal{U}}$ implies $u_{-1} \equiv e^{-w} \in \tilde{\mathcal{U}}$, we obtain the following equation:

$$\int Lw(x)\mu^*(dx) = 0, \qquad u = e^w \in \tilde{\mathcal{U}}.$$

Noting that $C_0^\infty(\mathbb{R}^N)$ is included in $\{w\colon u = e^w \in \tilde{\mathcal{U}}\}$, $\mu^*$ satisfies the following partial differential equation in distributional sense:

$$L^*\mu^* = 0 \qquad \text{in } \mathbb{R}^N,$$

where $L^*$ is a formal adjoint of $L$. Since we assumed the coefficients of $L$ are sufficiently smooth, $\mu^*$ has a density $p^*(x)$ and $p^*$ satisfies

$$L^*p^* = 0 \qquad \text{in } \mathbb{R}^N.$$

Here we recall that diffusion process $\{X_t\}$ does not explode in finite time because of Lemma 3.5. Then, by slight modifications of the Theorem in page 243 of [24] to the case that the second-order term of $L$ is divergence form, $\mu^*(dx) = p^*(x)\,dx$ is actually an invariant measure. □

Proof of Proposition 3.3. Let us define $U_0$ as follows:

$$U_0(x) = -(\tfrac{1}{2}a^{ij}D_{ij}W_0 + \tfrac{1}{2}\hat{a}\nabla W_0 \cdot \nabla W_0 + \tilde{b}\cdot\nabla W_0 + V),$$

where we take $W_0$ from (A3). By setting $W_{0,c} \equiv cW_0$ and $W_c \equiv cW$, we have

$$(3.12) \quad \frac{1}{2}a^{ij}D_{ij}W_{0,c} + \frac{1}{2c}\hat{a}\nabla W_{0,c} \cdot \nabla W_{0,c} + \tilde{b} \cdot \nabla W_{0,c} + cV = -cU_0,$$

$$\frac{1}{2}a^{ij}D_{ij}W_c + \frac{1}{2c}\hat{a}\nabla W_c \cdot \nabla W_c + \tilde{b} \cdot \nabla W_c + cV = c\Lambda,$$

where $c$ is in Remark 2.1. In the above equations, subtracting each side of the equations,

$$\frac{1}{2}a^{ij}D_{ij}(W_{0,c} - W_c) + (\tilde{b} + \hat{a}\nabla W) \cdot (\nabla W_{0,c} - \nabla W_c)$$

$$+ \frac{1}{2c}\hat{a}(\nabla W_{0,c} - \nabla W_c) \cdot (\nabla W_{0,c} - \nabla W_c) = -c(U_0 + \Lambda).$$

Define $\bar{\phi}$ as $\bar{\phi} = e^{W_{0,c} - W_c}$. Then, we have

$$(3.13) \quad \frac{1}{2}a^{ij}D_{ij}\bar{\phi} + (\tilde{b} + \hat{a}\nabla W) \cdot \nabla\bar{\phi} + \frac{1}{2c}((\hat{a} - ca)\nabla\bar{\phi} \cdot \nabla\bar{\phi})\frac{1}{\bar{\phi}} = -c(U_0 + \Lambda)\bar{\phi}.$$

Let $\{X_t\}$ be a solution of (3.1). By the Itô formula and (3.13),

$$d(\bar{\phi}(X_t)e^{\int_0^t c(U_0(X_s) + \Lambda)\,ds})$$

$$= \left[\frac{1}{2}a^{ij}D_{ij}\bar{\phi} + (\tilde{b} + \hat{a}\nabla W) \cdot \nabla\bar{\phi} + c(U_0 + \Lambda)\bar{\phi}\right](X_t)e^{\int_0^t c(U_0(X_s) + \Lambda)\,ds}\,dt$$

$$+ \sigma\nabla\bar{\phi}(X_t)e^{\int_0^t c(U_0(X_s) + \Lambda)\,ds}\,dB_t$$

$$= -\frac{1}{2c}\left[\frac{1}{\bar{\phi}}(\hat{a} - ca)\nabla\bar{\phi} \cdot \nabla\bar{\phi}\right](X_t)e^{\int_0^t c(U_0(X_s) + \Lambda)\,ds}\,dt$$

$$+ \sigma\nabla\bar{\phi}(X_t)e^{\int_0^t c(U_0(X_s) + \Lambda)\,ds}\,dB_t.$$

Since $ca(x) \le \hat{a}(x)$ and $\bar{\phi} > 0$, we obtain

$$(3.14) \quad \bar{\phi}(X_t)e^{\int_0^t c(U_0(X_s) + \Lambda)\,ds} \le \bar{\phi}(x) + \int_0^t \sigma\nabla\bar{\phi}(X_s)e^{\int_0^s c(U_0(X_r) + \Lambda)\,dr}\,dB_s.$$

By using stopping time $t \wedge \tau_n$ in (3.14),

$$E_x[\bar{\phi}(X_{t \wedge \tau_n})e^{\int_0^{t \wedge \tau_n} c(U_0(X_s) + \Lambda)\,ds}] \le \bar{\phi}(x).$$

Then, as $n \to \infty$, we have

$$(3.15) \quad E_x[\bar{\phi}(X_t)e^{\int_0^t c(U_0(X_s) + \Lambda)\,ds}; t < \tau_\infty] \le \bar{\phi}(x).$$

Let $\mathcal{C}_m$ be a subset in $\mathcal{M}_1(\mathbb{R}^N)$ defined as follows:

$$\mathcal{C}_m \equiv \{\mu \in \mathcal{M}_1(\mathbb{R}^N) : \mu(B_l) \ge 1 - \delta_l \ \forall\, l \ge m\}, \qquad m \ge 1,$$

where $\{\delta_l\}$ is a sequence such that $\delta_l \to 0$ and determined later. Note that $\mathcal{C}_m$ is a relative compact set in $\mathcal{M}_1(\mathbb{R}^N)$ because $\mathcal{C}_m$ is tight. From the definition of $T_t^{W,\Lambda}f$,

$$
\begin{aligned}
T_t^{W,\Lambda}f(x) &= E_x[f(X_t); \mu_t \in \mathcal{C}_m, t < \tau_\infty] \\
&\quad + E_x[f(X_t); \mu_t \notin \mathcal{C}_m, t < \tau_\infty] \\
&\leq \|f\|_\infty P_x[\mu_t \in \mathcal{C}_m, t < \tau_\infty] \\
&\quad + E_x[f(X_t); \mu_t \notin \mathcal{C}_m, t < \tau_\infty], \\
&\leq \|f\|_\infty P_x[\mu_t \in \mathcal{C}_m, t < \tau_\infty] \\
&\quad + \|f\bar{\phi}^{-1}\|_\infty E_x[\bar{\phi}(X_t); \mu_t \notin \mathcal{C}_m, t < \tau_\infty], \\
&\hspace{6cm} f \in C_0(\mathbb{R}^N), f \geq 0.
\end{aligned}
\tag{3.16}
$$

We shall prove that $E_x[\bar{\phi}(X_t); \mu_t \notin \mathcal{C}_m, t < \tau_\infty]$ decays exponentially as $t \to \infty$. On the event $\{\mu_t \notin \mathcal{C}_m, t < \tau_\infty\}$, there exists $l \geq m$ such that

$$
\mu_t(B_l) = \frac{1}{t}\int_0^t \mathbb{1}_{B_l}(X_s)\,ds \leq 1 - \delta_l
\tag{3.17}
$$

which is equivalent to

$$
\mu_t(B_l^c) = \frac{1}{t}\int_0^t \mathbb{1}_{B_l^c}(X_s)\,ds > \delta_l.
\tag{3.18}
$$

Then, we have on $\{\mu_t \notin \mathcal{C}_m, t < \tau_\infty\}$

$$
\begin{aligned}
\int_0^t c(U_0(X_s)+\Lambda)\,ds &= \int_0^t c(U_0(X_s)+\Lambda)\mathbb{1}_{B_l}(X_s)\,ds \\
&\quad + \int_0^t c(U_0(X_s)+\Lambda)\mathbb{1}_{B_l^c}(X_s)\,ds \\
&\geq \inf_x c(U_0(x)+\Lambda)\int_0^t \mathbb{1}_{B_l}(X_s)\,ds \\
&\quad + \inf_{|x|\geq l} c(U_0(x)+\Lambda)\int_0^t \mathbb{1}_{B_l^c}(X_s)\,ds \\
&= \beta_0 \mu_t(B_l)t + \beta_l \mu_t(B_l^c)t,
\end{aligned}
\tag{3.19}
$$

where we set $\beta_0 = \inf_x c(U_0(x)+\Lambda)$, $\beta_l = \inf_{|x|\geq l} c(U_0(x)+\Lambda)$. By (A3), there exists $m \geq 1$ such that

$$
\beta_l > 0 \qquad \forall l \geq m.
\tag{3.20}
$$

So, we obtain from (3.17)–(3.19),

$$
\int_0^t c(U_0(X_s)+\Lambda)\,ds \geq (-|\beta_0|(1-\delta_l)+\beta_l\delta_l)t.
$$

Take a positive constant $M > 0$. Then we choose $\delta_l$ such that $M = -|\beta_0|(1 - \delta_l) + \beta_l \delta_l$. Indeed, $\delta_l$ is defined by

$$\delta_l \equiv \frac{M + |\beta_0|}{|\beta_0| + \beta_l}.$$

Then, we have

$$(3.21) \qquad \int_0^t c(U_0(X_s) + \Lambda)\,ds \geq Mt \qquad \text{on } \{\mu_t \notin \mathcal{C}_m, t < \tau_\infty\}.$$

By (3.15) and (3.21),

$$\bar{\phi}(x) \geq E_x[\bar{\phi}(X_t)e^{\int_0^t c(U_0(X_s) + \Lambda)\,ds}; \mu_t \notin \mathcal{C}_m, t < \tau_\infty]$$
$$\geq e^{Mt}E_x[\bar{\phi}(X_t); \mu_t \notin \mathcal{C}_m, t < \tau_\infty].$$

Therefore we obtain

$$E_x[\bar{\phi}(X_t); \mu_t \notin \mathcal{C}_m, t < \tau_\infty] \leq \bar{\phi}(x)e^{-Mt}, \qquad t > 0.$$

By (3.16), we have

$$(3.22) \qquad T_t^{W,\Lambda}f(x) \leq \|f\|_\infty P_x[\mu_t \in \mathcal{C}_m, t < \tau_\infty] + \|f\bar{\phi}^{-1}\|_\infty \bar{\phi}(x)e^{-Mt}.$$

Applying Lemma 3.4,

$$\limsup_{t \to \infty} \frac{1}{t} \log P_x[\mu_t \in \mathcal{C}_m, t < \tau_\infty] \leq -\inf_{\mu \in \bar{\mathcal{C}}_m} I^W(\mu),$$

where $\bar{\mathcal{C}}_m$ is the closure of $\mathcal{C}_m$. Since $I^W(\mu)$ is lower semi-continuous and $\bar{\mathcal{C}}_m$ is compact in $\mathcal{M}_1(\mathbb{R}^N)$, $\inf_{\mu \in \bar{\mathcal{C}}_m} I^W(\mu)$ is attained at some $\mu^* \in \bar{\mathcal{C}}_m$. Since existence of invariant measure implies recurrence, it follows from Lemmas 3.5 and 3.6 and transience of $\{X_t\}$ that

$$\inf_{\mu \in \bar{\mathcal{C}}_m} I^W(\mu) > 0.$$

Hence, we can find a positive constant $\alpha_m > 0$ such that

$$(3.23) \qquad P_x[\mu_t \in \mathcal{C}_m, t < \tau_\infty] \leq C(x)e^{-\alpha_m t}, \qquad t > 0.$$

Then, from (3.22) and (3.23), we obtain

$$T_t^{W,\Lambda}f(x) \leq C(x)\|f\|_\infty e^{-\alpha_m t} + \|f\bar{\phi}^{-1}\|_\infty \bar{\phi}(x)e^{-Mt}. \qquad \square$$

We are ready to prove that for $\Lambda = \Lambda^*$, the corresponding diffusion process $\{X_t^*\}$ is ergodic.

THEOREM 3.7. *Let $(W^*, \Lambda^*)$ be a solution of (1.1) corresponding to $\Lambda^* = \inf \mathcal{A}$ and let $\{X_t^*\}$ be a solution of (3.1) for $(W^*, \Lambda^*)$. Under* (A1)–(A3), *$\{X_t^*\}$ is ergodic.*

PROOF.   Suppose that $\{X_t^*\}$ is transient. Then, by Proposition 3.3,

$$T_t^{W^*,\Lambda^*} f(x) \le C(x) e^{-\alpha t} \qquad \forall f \in C_0(\mathbb{R}^N), f \ge 0.$$

Note that $\alpha$ is a positive constant independent of $f$ and $x$. Taking $0 < \varepsilon < \alpha$, we see that

$$\int_0^\infty E_x[f(X_t^*)e^{\varepsilon t}; t < \tau_\infty]\,dt = \int_0^\infty T_t^{W^*,\Lambda^*} f(x)e^{\varepsilon t}\,dt$$

$$= C(x)\int_0^\infty e^{-(\alpha-\varepsilon)t}\,dt < \infty.$$

Then, there exists Green function $G(x,y)$ for $(1/2)a^{ij}D_{ij} + (\tilde{b} + \hat{a}\nabla W^*) \cdot \nabla + \varepsilon$ and $G(x,y)$ satisfies the following:

$$(3.24)\quad \tfrac{1}{2}a^{ij}D_{ij}G(\cdot,y) + (\tilde{b} + \hat{a}\nabla W^*) \cdot \nabla G(\cdot,y) + \varepsilon G(\cdot,y) = 0 \qquad \text{in } \mathbb{R}\backslash\{y\}.$$

We take a sequence $\{y_n\}$ in $\mathbb{R}^N$ such that $y_n \in B_{n+1}\backslash \bar{B}_n$. Define $\bar\phi_n(x)$ as follows:

$$\bar\phi_n(x) \equiv \frac{G(x,y_n)}{G(0,y_n)}, \qquad x \in \mathbb{R}^N\backslash\{y_n\}.$$

Then, we have from (3.24)

$$(3.25)\quad \tfrac{1}{2}a^{ij}D_{ij}\bar\phi_n + (\tilde{b} + \hat{a}\nabla W^*) \cdot \nabla\bar\phi_n + \varepsilon\bar\phi_n = 0 \qquad \text{in } \mathbb{R}^N\backslash\{y_n\}.$$

We note that by setting $\bar W_n \equiv (1/\bar c)\log\bar\phi_n$, (3.25) is equivalent to the following:

$$\tfrac{1}{2}a^{ij}D_{ij}\bar W_n + \frac{\bar c}{2}a^{ij}D_i\bar W_n D_j\bar W_n + (\tilde{b} + \hat{a}\nabla W^*) \cdot \nabla\bar W_n + \frac{\varepsilon}{\bar c} = 0 \qquad \text{in } \mathbb{R}^N\backslash\{y_n\},$$

where $\bar c$ is taken from Remark 2.1. By Lemma 2.4, we have

$$\sup_{B_r}|\nabla\bar W_n| \le C_r, \qquad r < n.$$

Thus, in a similar way to the proof of existence of solutions of (1.1), we can see that there exists smooth function $\bar W$ such that

$$(3.26)\quad \tfrac{1}{2}a^{ij}D_{ij}\bar W + (\tilde{b} + \hat{a}\nabla W^*) \cdot \nabla\bar W + \frac{\bar c}{2}a\nabla\bar W \cdot \nabla\bar W + \frac{\varepsilon}{\bar c} = 0.$$

Since $(W^*,\Lambda^*)$ is a solution of (1.2),

$$(3.27)\quad \tfrac{1}{2}a^{ij}D_{ij}W^* + \tilde{b}\cdot\nabla W^* + \tfrac{1}{2}\hat{a}\nabla W^* \cdot \nabla W^* + V - \Lambda^* = 0.$$

Adding (3.26) to (3.27), it follows from Remark 2.1 that

$$0 = \frac{1}{2}a^{ij}D_{ij}(W^* + \bar W) + \tilde{b}\cdot(\nabla W^* + \nabla\bar W)$$

$$+ \frac{1}{2}\hat{a}\nabla W^* \cdot \nabla W^* + \hat{a}\nabla W^* \cdot \nabla \bar{W} + \frac{\bar{c}}{2}a\nabla\bar{W} \cdot \nabla\bar{W} + V - \left(\Lambda^* - \frac{\varepsilon}{\bar{c}}\right)$$

$$\geq \frac{1}{2}a^{ij}D_{ij}(W^* + \bar{W}) + \tilde{b}\cdot(\nabla W^* + \nabla\bar{W})$$

$$+ \frac{1}{2}\hat{a}\nabla W^* \cdot \nabla W^* + \hat{a}\nabla W^* \cdot \nabla\bar{W} + \frac{1}{2}\hat{a}\nabla\bar{W} \cdot \nabla\bar{W} + V - \left(\Lambda^* - \frac{\varepsilon}{\bar{c}}\right)$$

$$= \frac{1}{2}a^{ij}D_{ij}(W^* + \bar{W}) + \tilde{b}\cdot\nabla(W^* + \bar{W})$$

$$+ \frac{1}{2}\hat{a}\nabla(W^* + \bar{W}) \cdot \nabla(W^* + \bar{W}) + V - \left(\Lambda^* - \frac{\varepsilon}{\bar{c}}\right).$$

Thus, $W^* + \bar{W}$ is a super solution of (1.1) for $\Lambda = \Lambda^* - \varepsilon/\bar{c}$. In the same way as the proof that $\mathcal{A} \neq \varnothing$ given in Section 2 we can see that there exists a smooth function $\tilde{W}$ such that

$$\frac{1}{2}a^{ij}D_{ij}\tilde{W} + \frac{1}{2}\hat{a}\nabla\tilde{W} \cdot \nabla\tilde{W} + \tilde{b}\cdot\nabla\tilde{W} + V = \Lambda^* - \frac{\varepsilon}{\bar{c}}.$$

This contradicts $\Lambda^* = \inf\mathcal{A}$. Therefore, $\{X_t^*\}$ is recurrent.

In order to see that $\{X_t^*\}$ is actually ergodic, we recall the proof of Proposition 3.3. If we suppose $\inf_{\mu\in\bar{\mathcal{C}}_m} I^{W^*}(\mu) > 0$ where $m$ is chosen in (3.20), we can prove Proposition 3.3, which implies that $\{X_t^*\}$ is transient. Hence, we see that

$$(3.28) \qquad\qquad \inf_{\mu\in\bar{\mathcal{C}}_m} I^{W^*}(\mu) = 0.$$

Since $\bar{\mathcal{C}}_m$ is compact, (3.28) is attained at $\mu^* \in \bar{\mathcal{C}}_m$. Then, it follows from Lemmas 3.5 and 3.6 that $\mu^*$ is an invariant measure for $\{X_t^*\}$. □

3.2. *Uniqueness of solutions corresponding to the bottom.* We proved that solution $(W^*, \Lambda^*)$ of (1.1) for $\Lambda^* = \inf\mathcal{A}$ corresponds to ergodicity to $\{X_t^*\}$ of (3.1). Now we shall show that the solution corresponding to $\Lambda^*$ is unique up to an additive constant. Note that the solution of (1.1) has ambiguity on an additive constant.

THEOREM 3.8. *Let $W_i^*$, $i = 1, 2$, be solutions of (1.1) corresponding to $\Lambda^* = \inf\mathcal{A}$. Under* (A1)–(A3), *there exists constant $k$ such that $W_2^*(x) = W_1^*(x) + k$.*

PROOF. Since $W_i^*$, $i = 1, 2$, are solutions of (1.2),

$$\frac{1}{2}a^{ij}D_{ij}W_1^* + \frac{1}{2}\hat{a}\nabla W_1^* \cdot \nabla W_1^* + \tilde{b}\cdot\nabla W_1^* + V = \Lambda^*,$$

$$\frac{1}{2}a^{ij}D_{ij}W_2^* + \frac{1}{2}\hat{a}\nabla W_2^* \cdot \nabla W_2^* + \tilde{b}\cdot\nabla W_2^* + V = \Lambda^*.$$

Subtracting each side in the above equations, we have

$$
\begin{aligned}
(3.29) \quad & \tfrac{1}{2} a^{ij} D_{ij}(W_1^* - W_2^*) + (\tilde{b} + \hat{a}\nabla W_2^*) \cdot (\nabla W_1^* - \nabla W_2^*) \\
& + \tfrac{1}{2}\hat{a}(\nabla W_1^* - \nabla W_2^*) \cdot (\nabla W_1^* - W_2^*) = 0.
\end{aligned}
$$

Let us set $\phi(x) \equiv e^{c(W_1^*(x) - W_2^*(x))}$, where $c$ is in Remark 2.1. Rewriting (3.29) in terms of $\phi$, we have

$$
\tfrac{1}{2} a^{ij} D_{ij}\phi + (\tilde{b} + \hat{a}\nabla W_2^*) \cdot \nabla\phi + \frac{1}{2c}(\hat{a} - ca)\frac{\nabla\phi}{\phi} \cdot \nabla\phi = 0.
$$

Hence it is implied from Remark 2.1 that

$$
(3.30) \qquad L\phi \equiv \tfrac{1}{2} a^{ij} D_{ij}\phi + (\tilde{b} + \hat{a}\nabla W_2^*) \cdot \nabla\phi \le 0.
$$

Let us take $x, y \in \mathbb{R}^N$ and consider the SDE of (3.1) for $W = W_2^*$:

$$
dX_t^* = (\tilde{b}(X_t^*) + \hat{a}\nabla W_2^*(X_t^*))\, dt + \sigma(X_t^*)\, dB_t, \qquad X_0^* = x.
$$

Define $\tau_{B_n} = \inf\{t : X_t^* \notin B_n\}$, $\sigma_{B_\varepsilon(y)} = \inf\{t : X_t^* \in B_\varepsilon(y)\}$. Note that $\{X_t^*\}$ is ergodic from Theorem 3.7, especially recurrent. It follows from the Itô formula and (3.30) that

$$
\begin{aligned}
\phi(X_{t \wedge \tau_{B_n} \wedge \sigma_{B_\varepsilon(y)}}^*) &= \phi(x) + \int_0^{t \wedge \tau_{B_n} \wedge \sigma_{B_\varepsilon(y)}} L\phi(X_s^*)\, ds \\
& \quad + \int_0^{t \wedge \tau_{B_n} \wedge \sigma_{B_\varepsilon(y)}} \nabla\phi(X_s^*) \cdot \sigma(X_s^*)\, dB_s \\
& \le \phi(x) + \int_0^{t \wedge \tau_{B_n} \wedge \sigma_{B_\varepsilon(y)}} \nabla\phi(X_s^*) \cdot \sigma(X_s^*)\, dB_s.
\end{aligned}
$$

Thus we have $E_x[\phi(X_{t \wedge \tau_{B_n} \wedge \sigma_{B_\varepsilon(y)}}^*)] \le \phi(x)$. By taking the limit as $n \to \infty$, it follows by Fatou's lemma that $E[\phi(X_{t \wedge \sigma_{B_\varepsilon(y)}}^*)] \le \phi(x)$. Noting that $P_x[\sigma_{B_\varepsilon(y)} < \infty] = 1$, we have by sending $t \to \infty$,

$$
E_x[\phi(X_{\sigma_{B_\varepsilon(y)}}^*)] \le \phi(x).
$$

We note again that $\{X_t^*\}$ hits the boundary of $B_\varepsilon(y)$ in finite time with probability 1. Hence we can see that

$$
\phi(x) \ge E_x[\phi(X_{\sigma_{B_\varepsilon(y)}}^*)] \ge \inf_{\partial B_\varepsilon(y)} \phi.
$$

Taking the limit as $\varepsilon \to 0$, we obtain

$$
\phi(y) \le \phi(x), \qquad x, y \in \mathbb{R}^N,
$$

which implies $\phi$ is constant. Therefore $W_1^* - W_2^*$ is also constant. $\quad\square$

EXAMPLE 3.9. Let us consider the linear case:

$$b(x) = Dx;$$
$$a(x) = a, \qquad \hat{a}(x) = \hat{a};$$
$$V(x) = \tfrac{1}{2}x \cdot Mx + v \cdot x.$$

$D, M, a, \hat{a}$ are matrices and $M$ is symmetric; $a, \hat{a}$ are positive-definite. We consider quadratic solution $W$,

$$W(x) = \tfrac{1}{2}x \cdot Kx + e \cdot x.$$

Then

$$K\hat{a}K + D^T K + KD + M = 0,$$
$$(D^T + K\hat{a})e + v = 0,$$
$$\Lambda = \tfrac{1}{2}\operatorname{tr}(aK) + \tfrac{1}{2}e \cdot \hat{a}e.$$

If $M$ is negative-definite, then there is a unique solution $K$ such that $K$ is nonpositive-definite and $D + \hat{a}K$ is stable. See [15]. For such $K$, the equation for $e$ can be uniquely solved. The stability of $D + \hat{a}K$ implies that the diffusion

$$dX_t = (D + \hat{a}K)X_t \, dt + \sigma \, dB_t$$

is ergodic, $\sigma$ is the square root of $a$. This implies $\Lambda$ defined above is equal to $\Lambda^*$ and $W$ obtained is the unique solution corresponding to $\Lambda^*$. We note there are also solutions that are not quadratic.

In particular, in the one-dimensional case, the equation for $K$ becomes

$$\hat{a}K^2 + 2DK + M = 0.$$

If $M < 0$, the solution is given by

$$K = -\frac{D}{\hat{a}} \pm \sqrt{\left(\frac{D}{\hat{a}}\right)^2 - M}.$$

Let

$$K_- = -\frac{D}{\hat{a}} - \sqrt{\left(\frac{D}{\hat{a}}\right)^2 - M},$$

$$K_+ = -\frac{D}{\hat{a}} + \sqrt{\left(\frac{D}{\hat{a}}\right)^2 - M}.$$

Then

$$D + \hat{a}K_- = -\sqrt{\left(\frac{D}{\hat{a}}\right)^2 - M} < 0$$

and

$$D + \hat{a}K_+ = \sqrt{\left(\frac{D}{\hat{a}}\right)^2 - M} > 0.$$

Solution corresponding to $K_-$ is $W^*$.

EXAMPLE 3.10.  In [9, 10], the following conditions are considered:

(a) $a(x) = \hat{a}(x) = I$.

(b) $b(0) = 0$; $b_j(x)$ has continuous first-order derivatives, $D_i b_j(x)$ is bounded for all $i, j$.

(c) There is $c_0 > 0$ such that for all $x, \eta \in \mathbb{R}^d$, $\eta \cdot Db(x)\eta \leq -c_0|\eta|^2$. Here $Db(x) = (D_i b_j(x))_{ij}$.

(d) $V(x), D_i V(x), i = 1, \ldots, d$, are bounded.

Under these conditions, Fleming and co-workers prove that there exists unique solution $(W, \Lambda)$ satisfying the condition

$$|\nabla W(x)| \leq \frac{1}{c_0}\|\nabla V\|_\infty,$$

$\|\nabla V\|_\infty = \sup_x |\nabla V(x)|$. Let $\{X_t\}$ be the diffusion,

$$dX_t = (b(X_t) + \nabla W(X_t)) \, dt + dB_t.$$

Denote

$$c = \frac{1}{c_0}\|\nabla V\|_\infty.$$

Then

$$d|X_t|^2 = 2\left(X_t \cdot b(X_t) + X_t \cdot \nabla W(X_t) + \frac{d}{2}\right) dt + 2X_t \, dB_t.$$

By (a), (c) and the mean value theorem,

$$x \cdot b(x) \leq -c_0|x|^2.$$

Property of $W$ implies

$$x \cdot \nabla W(x) \leq c|x|.$$

By using a routine argument and considering $|X_t|^2 \exp(c_0 t)$, we can prove

$$E_x[|X_t|^2] \leq \exp(-c_0 t)|x|^2 + \frac{1}{c_0}\left(d + \frac{c^2}{c_0}\right).$$

This implies $\{X_t\}$ is ergodic. Therefore, $\Lambda = \Lambda^*$ and $W = W^*$.

In [15], different conditions are considered that are given as follows:

(a)$'$ There are $c_1, c_2 > 0$ such that

$$c_1 \leq a(x) \leq c_2, \qquad c_1 \leq \hat{a}(x) \leq c_2.$$

(b)$'$ There are $c_0, r_0 > 0$ such that

$$x \cdot b(x) \leq -c_0 |x|^2, \qquad |x| \geq r_0.$$

(c)$'$ $a_{ij}(x), \hat{a}_{ij}(x)$ and $b_i(x)$ have bounded first-order derivatives.

(d)$'$ $V(x), \nabla V(x)$ are bounded.

Then (1.1) with $\Lambda = \Lambda^*$ has unique solution $W^*$ with $W^*(0) = 0$. Moreover, for any $\alpha, \beta > 0$, there are $c_{\alpha, \beta}$ such that

$$|W^*(x)| \leq \alpha |x|^\beta + c_{\alpha, \beta},$$

$$|\nabla W^*(x)| \leq \alpha |x|^\beta + c_{\alpha, \beta}.$$

## 4. Perturbation of coefficients.

In the present section we shall consider the structures of solutions of (1.1) under perturbation of coefficients. This is to consider the following equation parameterized by $n \in \mathbb{N}$:

$$(4.1) \quad \tfrac{1}{2} D_i(a_n^{ij} D_j W_n) + \tfrac{1}{2} \hat{a}_n^{ij} D_i W_n D_j W_n + b_n \cdot \nabla W_n + V_n = \Lambda_n \qquad \text{in } \mathbb{R}^N,$$

or equivalently

$$\tfrac{1}{2} a_n^{ij} D_{ij} W_n + \tfrac{1}{2} \hat{a}_n^{ij} D_i W_n D_j W_n + \tilde{b}_n \cdot \nabla W_n + V_n = \Lambda_n,$$
$$(4.2)$$
$$\tilde{b}_n^i(x) \equiv b_n^i(x) + \tfrac{1}{2} D_j \hat{a}_n^{ij}(x).$$

In the same way as (2.1) we define $\mathcal{A}_n$ as follows:

$$(4.3) \quad \mathcal{A}_n \equiv \{\Lambda_n; \text{ there exists smooth } W_n \text{ satisfying (4.1) for } \Lambda_n\}.$$

Under certain assumptions, we can show that $\mathcal{A}_n = [\Lambda_n^*, \infty)$ for some finite $\Lambda_n^*$. Furthermore, we can classify the set $\mathcal{A}_n$ of $\Lambda_n$ according to the global behavior of diffusion process defined by the SDE:

$$(4.4) \quad dX_t = (\tilde{b}_n(X_t) + \hat{a}_n \nabla W_n(X_t)) \, dt + \sigma_n(X_t) \, dB_t, \qquad \sigma_n(x) \equiv a_n(x)^{1/2},$$

where $W_n$ is a solution of (4.1) corresponding to $\Lambda_n$ and $\{B_t\}$ is $\mathcal{F}_t$-standard Brownian motion on a filtered probability space $(\Omega, \mathcal{F}, P, \{\mathcal{F}_t\})$. Indeed, we can prove that $\Lambda_n = \Lambda_n^*$ (resp. $\Lambda_n > \Lambda_n^*$) corresponds to ergodicity (resp. transience) of $\{X_t\}$ defined by (4.4) and $W_n^*$ corresponding to $\Lambda_n^*$ is unique up to additive constants.

It is interesting to study stability of solution $(W_n^*, \Lambda_n^*)$ corresponding to the bottom of $\mathcal{A}_n$ under perturbation of coefficients. Suppose that all coefficients converge to corresponding ones, respectively, in some sense:

$$a_n^{ij} \to a^{ij}, \qquad \hat{a}_n^{ij} \to \hat{a}^{ij}, \qquad b_n \to b, \qquad V_n \to V \qquad \text{as } n \to \infty.$$

We hope to prove that $W_n^*$ and $\Lambda_n^*$ converge to $W^*$ and $\Lambda^* = \inf \mathcal{A}$, respectively, where $\mathcal{A}$ is defined by (2.1) and $W^*$ is a unique solution of (1.1) corresponding to $\Lambda^*$. This means that the solutions corresponding to $\Lambda_n^* = \inf \mathcal{A}_n$ are stable under the perturbation.

It turns out this is a delicate problem. For the illustration of the idea, we shall be content with the following special example. The result obtained will be used in Section 5. We refer to [16] for more general discussion and a counterexample.

We now consider the following special example; we consider the following equation:

$$\tfrac{1}{2}a^{ij}D_{ij}W_n + \tfrac{1}{2}\hat{a}^{ij}D_iW_nD_jW_n + \tilde{b}\cdot\nabla W_n + V_n = \Lambda_n, \qquad V_n = V_0 + \bar{V}_n.$$
(4.5)

The equation corresponding to the limit of (4.5) is as follows:

$$(4.6) \qquad \tfrac{1}{2}a^{ij}D_{ij}W + \tfrac{1}{2}\hat{a}^{ij}D_iWD_jW + \tilde{b}\cdot\nabla W + V = \Lambda, \qquad V = V_0 + \bar{V}.$$

We assume the following conditions. We suppose implicitly that all the coefficients are smooth.

(B1) There exists smooth function $W_0$ such that

$$\bar{U}_0 \equiv -(\tfrac{1}{2}a^{ij}D_{ij}W_0 + \tfrac{1}{2}\hat{a}^{ij}D_iW_0D_jW_0 + \tilde{b}\cdot\nabla W_0 + V_0) \to \infty$$
$$\text{as } |x| \to \infty.$$

(B2) $\bar{V}_n$, $D\bar{V}_n$ are bounded in $\mathbb{R}^N$ uniformly on $n$.

(B3) $\bar{V}_n$ converges to $\bar{V}$ uniformly on each compact set as $n \to \infty$.

Let $W_n^*$ (resp. $W^*$) be solution of (4.5) [resp. (4.6)] corresponding to $\Lambda_n^*$ (resp. $\Lambda^*$).

PROPOSITION 4.1.   *Under* (A1), (A2), (B1)–(B3), $\Lambda_n^*$ *converges to* $\Lambda^*$.

PROOF.   It is easy to prove that $\liminf_{n\to\infty} \Lambda_n^* \geq \Lambda^*$. We shall prove $\limsup_{n\to\infty} \Lambda_n^* \leq \Lambda^*$. Since $W_n^*$ (resp. $W^*$) is solution of (4.5) for $\Lambda_n^*$ [resp. (4.6) for $\Lambda^*$],

$$\tfrac{1}{2}a^{ij}D_{ij}W_n^* + \tfrac{1}{2}\hat{a}^{ij}D_iW_n^*D_jW_n^* + \tilde{b}\cdot\nabla W_n^* + V_0 + \bar{V}_n = \Lambda_n^*,$$
$$\tfrac{1}{2}a^{ij}D_{ij}W^* + \tfrac{1}{2}\hat{a}^{ij}D_iW^*D_jW^* + \tilde{b}\cdot\nabla W^* + V_0 + \bar{V} = \Lambda^*.$$

By subtracting both of sides, we have

$$\tfrac{1}{2}a^{ij}D_{ij}(W^* - W_n^*) + (\tilde{b} + \hat{a}\nabla W_n^*)\cdot\nabla(W^* - W_n^*)$$
$$+ \tfrac{1}{2}\hat{a}\nabla(W^* - W_n^*)\cdot\nabla(W^* - W_n^*) = \Lambda^* - \Lambda_n^* + \bar{V}_n - \bar{V}.$$

Then, we can see that for $c > 0$,

$$\frac{1}{2}a^{ij}D_{ij}(c(W^* - W_n^*)) + (\tilde{b} + \hat{a}\nabla W_n^*) \cdot \nabla(c(W^* - W_n^*))$$

$$+ \frac{1}{2}a\nabla(c(W^* - W_n^*)) \cdot \nabla(c(W^* - W_n^*))$$

$$= c(\Lambda^* - \Lambda_n^* + \bar{V}_n - \bar{V}) + \frac{1}{2c}(ca - \hat{a})\nabla(c(W^* - W_n^*)) \cdot \nabla(c(W^* - W_n^*)).$$

If we take $c > 0$ from Remark 2.1, then

$$(4.7) \quad \begin{aligned} &\frac{1}{2}a^{ij}D_{ij}(c(W^* - W_n^*)) + (\tilde{b} + \hat{a}\nabla W_n^*) \cdot \nabla(c(W^* - W_n^*)) \\ &+ \frac{1}{2}a\nabla(c(W^* - W_n^*)) \cdot \nabla(c(W^* - W_n^*)) \le c(\Lambda^* - \Lambda_n^* + \bar{V}_n - \bar{V}). \end{aligned}$$

Let $\mu_n^*$ be invariant measure of diffusion process $X_t^{*,n}$ defined by the following SDE:

$$dX_t^{*,n} = (b + \hat{a}\nabla W_n^*)(X_t^{*,n})\,dt + \sigma(X_t^{*,n})\,dB_t, \qquad X_0^{*,n} = x.$$

By the construction of $\mu_n^*$, we see that

$$(4.8) \qquad\qquad I_n^{W_n^*}(\mu_n^*) = 0,$$

where

$$I_n^{W_n^*}(\mu) = -\inf_{u \in \mathcal{U}_n} \int \frac{L_n u}{u}\,d\mu, \qquad \mu \in \mathcal{M}_1,$$

$$L_n = \frac{1}{2}a^{ij}D_{ij} + (b + \hat{a}\nabla W_n^*) \cdot \nabla,$$

$$\mathcal{U}_n = \{u \in C^2(\mathbb{R}^N); u(x) > 0, L_n u/u \text{ is bounded above}\}.$$

Set $u = \exp(c(W^* - W_n^*))$. Then, from (4.7), we have

$$(4.9) \quad \begin{aligned} \frac{L_n u}{u} &= \frac{1}{2}a^{ij}D_{ij}(c(W^* - W_n^*)) + (\tilde{b} + \hat{a}\nabla W_n^*) \cdot \nabla(c(W^* - W_n^*)) \\ &\quad + \frac{1}{2}a\nabla(c(W^* - W_n^*)) \cdot \nabla(c(W^* - W_n^*)) \\ &\le c(\Lambda^* - \Lambda_n^* + \bar{V}_n - \bar{V}). \end{aligned}$$

Since $\bar{V}$, $\bar{V}_n$ are bounded, $u \in \mathcal{U}_n$. Hence, by (4.8),

$$\int \frac{L_n u}{u}\,d\mu_n^* \ge 0.$$

Then, the above inequality and (4.9) imply that

$$c \int (\Lambda^* - \Lambda_n^* + \bar{V}_n - \bar{V})\,d\mu_n^* \ge 0.$$

Therefore we have

$$\Lambda^* - \Lambda_n^* \geq \int (\bar{V} - \bar{V}_n)\, d\mu_n^*.$$

We note that $\{\mu_n^*\}$ is tight by the proof of Proposition 3.3. Indeed, $\{\mu_n^*\} \in \bar{\mathcal{C}}_m$,

$$\mathcal{C}_m \equiv \{\mu \in \mathcal{M}_1(\mathbb{R}^N); \mu(B_l) \geq 1 - \delta_l \ \forall\, l \geq m\}, \qquad \delta_l \equiv \frac{M + |\beta_0|}{|\beta_0| + \beta_l},$$

$$\beta_0 \equiv \inf_{x \in \mathbb{R}^N} \{U_0(x) + \alpha\}, \qquad \beta_l \equiv \inf_{|x| \geq l} \{U_0(x) + \alpha\},$$

$$U_0(x) \equiv -\sup_n (\tfrac{1}{2} a^{ij} D_{ij} W_0 + \tfrac{1}{2} \hat{a} \nabla W_0 \cdot \nabla W_0 + \tilde{b} \cdot \nabla W_0 + V_n),$$

$$\alpha = \inf_n \Lambda_n^*.$$

$W_0$ is from (B1) and $m$ is taken so that

$$\beta_l > 0 \qquad \forall\, l \geq m.$$

Since $\bar{V}_n$ are uniformly bounded and $\bar{V}_n$ converges to $\bar{V}$ uniformly on compact sets, we can see that

$$\int (\bar{V} - \bar{V}_n)\, d\mu_n^* \to 0 \qquad (n \to \infty).$$

Therefore, we have

$$\limsup_{n \to \infty} \Lambda_n^* \leq \Lambda^*. \qquad \qquad \square$$

THEOREM 4.2. *Let $W_n^*$ (resp. $W^*$) be solution of (4.5) for $\Lambda_n^*$ [resp. (4.6) for $\Lambda^*$]. Under (A1), (A2), (B1)–(B3), $W_n^*$ (resp. $\Lambda_n^*$) converges to $W^*$ (resp. $\Lambda^*$) uniformly on compact sets, $H^1_{\mathrm{loc}}$-strongly as $n \to \infty$.*

PROOF. By using the estimate of $\nabla W_n$ as in Lemma 2.4 and the argument in the proof of Theorem 2.6, we can show $W_n^*$ (resp. $\Lambda_n^*$) converges to $\tilde{W}^*$ (resp. $\tilde{\Lambda}^*$) uniformly on compact sets, $H^1_{\mathrm{loc}}$-strongly by taking a subsequence if necessary. Furthermore $(\tilde{W}^*, \tilde{\Lambda}^*)$ is a solution of (4.6). Indeed, we see that $\Lambda^* = \tilde{\Lambda}^*$ by Proposition 4.1. By uniqueness of solution corresponding to $\Lambda^*$ in Proposition 4.1, we have $W^*(x) = \tilde{W}^*(x)$.  $\square$

**5. Representation of $\Lambda^*$.** For (1.1) with the coefficients satisfying the conditions (A1)–(A3), we have proved that there is a unique solution $(W^*, \Lambda^*)$ with $W^*(0) = 0$ and $\Lambda = \Lambda^*$ is the smallest such that (1.1) has solution $W$. In this section we will give a representation of $\Lambda^*$. From the representation, we will get some moment condition for $\mu^*$, the invariant measure for the diffusion in (3.1) constructed from $W = W^*$. Before we state our main results, we give some notation.

We shall consider a family of $V$. That is, we consider a particular $V_0$ and the bounded perturbation of $V_0$, say $V_0 + \bar{V}$ for bounded $\bar{V}$. Therefore, in this section we consider the equation

$$(5.1) \quad \tfrac{1}{2} D_i(a^{ij} D_j W) + \tfrac{1}{2} \hat{a}^{ij} D_i W D_j W + b \cdot \nabla W + V_0 + V = \Lambda \qquad \text{in } \mathbb{R}^N.$$

The smallest $\Lambda$ such that this has solution $W$ is denoted $\Lambda^*(V)$. We still use $W^*$ for the solution corresponding to $\Lambda^*(V)$. We shall mention if we want to emphasize the dependence of $W^*$ on $V$.

In this section we assume the following condition:

(A3)$'$  $V_0$ is smooth and there exists a smooth function $W_0$ such that

$$U_0(x) = -(\tfrac{1}{2} D_i(a^{ij} D_j W_0) + b \cdot \nabla W_0 + \tfrac{1}{2} \hat{a}^{ij} D_i W_0 D_j W_0 + V_0) \to \infty$$
$$\text{as } |x| \to \infty.$$

We consider $V$ satisfying the condition:

$$(5.2) \qquad V \text{ is smooth}, \qquad |V(x)| \le c, \qquad |DV(x)| \le c,$$

where $c$ is a constant that may depend on $V$.

Let $W$ be a smooth function. We denote

$$(5.3) \qquad G(W) = \tfrac{1}{2} D_i(a^{ij} D_j W) + b \cdot \nabla W + \tfrac{1}{2} \hat{a}^{ij} D_i W D_j W + V_0.$$

For each probability measure $\mu$ on $\mathbb{R}^N$, we define

$$J(\mu) = \sup\{-\langle G(W), \mu \rangle; W \text{ is smooth and } G(W) \text{ is bounded above}\}.$$

Here for a function $f$

$$\langle f, \mu \rangle = \int f(x) \, d\mu(x).$$

Now we can state our main results.

THEOREM 5.1.  *Assume* (A1), (A2), (A3)$'$ *and* (5.2). *Let* $(W^*, \Lambda^*(V))$ *be the solution defined above. Then* $\Lambda^*(V)$ *has the representation*

$$\Lambda^*(V) = \sup_{\mu \in \mathcal{M}_1(\mathbb{R}^N)} \{\langle V, \mu \rangle - J(\mu)\}.$$

*Here* $\mathcal{M}_1(\mathbb{R}^N)$ *denotes the set of all probability measures on* $\mathbb{R}^N$. *The supremum attains at* $\mu = \mu^*$, $\mu^*$ *is the invariant measure of the diffusion in* (3.1) *for* $W = W^*$.

From this theorem, we have

$$\Lambda^*(V) = \langle V, \mu^* \rangle - J(\mu^*)$$
$$\leq \langle V, \mu^* \rangle + \langle G(W_0), \mu^* \rangle$$
$$= \langle V, \mu^* \rangle + \langle -U_0, \mu^* \rangle.$$

Therefore, we have the following corollary.

COROLLARY 5.2. *We assume the condition as in the above theorem. Then*

$$\int U_0(x) \, d\mu^*(x) \leq -\Lambda^*(V) + \|V\|_\infty.$$

*In particular,*

$$\int U_0(x) \, d\mu^*(x) < \infty.$$

*Here $\|V\|$ is the supnorm of $V$.*

Before we prove Theorem 5.1, we mention some elementary properties of $\Lambda^*(V)$.

LEMMA 5.3. (i) $\Lambda^*(V)$ *is Lipschitz with constant 1. That is,*

$$|\Lambda^*(V_1) - \Lambda^*(V_2)| \leq \|V_1 - V_2\|_\infty$$

*for $V_1, V_2$ satisfying (5.2). From this, $\Lambda^*(V)$ can be defined for all bounded continuous functions $V$ by extension.*

(ii) $\Lambda^*(V)$ *is convex in $V$.*

PROOF. Let $W_1^*$ be the solution of (5.1) for $V = V_1, \Lambda = \Lambda^*(V_1)$. Then

$$\tfrac{1}{2} D_i(a^{ij} D_j W_1^*) + \tfrac{1}{2} \hat{a}^{ij} D_i W_1^* D_j W_1^* + b \cdot \nabla W_1^* + V_0 + V_2$$
$$= \Lambda^*(V_1) + V_2 - V_1$$
$$\leq \Lambda^*(V_1) + \|V_2 - V_1\|_\infty.$$

This implies

$$\Lambda^*(V_2) \leq \Lambda^*(V_1) + \|V_2 - V_1\|_\infty.$$

See the argument in the proof of Theorem 2.6. Similarly,

$$\Lambda^*(V_1) \leq \Lambda^*(V_2) + \|V_2 - V_1\|.$$

Therefore,

$$|\Lambda^*(V_1) - \Lambda^*(V_2)| \leq \|V_1 - V_2\|_\infty.$$

We now show that $\Lambda^*(V)$ is convex. Let $V_1, V_2$ satisfy (5.2) and $W_k^*$, $k = 1, 2$, satisfy

$$\tfrac{1}{2}D_i(a^{ij}D_jW_k^*) + \tfrac{1}{2}\hat{a}^{ij}D_iW_k^*D_jW_k^* + b\cdot\nabla W_k^* + V_0 + V_k = \Lambda^*(V_k), \qquad k = 1, 2.$$

Let $0 < \lambda < 1$ and denote $W = \lambda W_1^* + (1 - \lambda)W_2^*$, $V = \lambda V_1 + (1 - \lambda)V_2$ and $\Lambda = \lambda\Lambda^*(V_1) + (1 - \lambda)\Lambda^*(V_2)$. Then by a simple calculation, we have

$$\tfrac{1}{2}D_i(a^{ij}D_jW) + \tfrac{1}{2}\hat{a}^{ij}D_iWD_jW + b\cdot\nabla W + V_0 + V \leq \Lambda, \qquad k = 1, 2.$$

This implies

$$\Lambda^*(V) \leq \Lambda.$$

See the proof of Theorem 2.6. That is, we have

$$\Lambda^*(\lambda V_1 + (1 - \lambda)V_2) \leq \lambda\Lambda^*(V_1) + (1 - \lambda)\Lambda^*(V_2). \qquad \square$$

Denote by $C_b(\mathbb{R}^N)$ the collection of all bounded continuous functions defined on $\mathbb{R}^N$. $C_b(\mathbb{R}^N)$ is a Banach space with supnorm. The dual space $C_b(\mathbb{R}^N)^*$ can be identified with the set of all regular bounded finitely additive set functions defined on the field generated by closed sets of $\mathbb{R}^N$. That is, for an element $T \in C_b(\mathbb{R}^N)^*$, there is regular bounded finitely additive set function $\mu$ such that

$$T(V) = \int V(x)\, d\mu(x), \qquad V \in C_b(\mathbb{R}^n).$$

See [6], Theorem IV.6.2. For regular additive function see [6], Theorem III.5.11. We note that $\mathcal{M}_1(\mathbb{R}^N)$ is a subset of $C_b(\mathbb{R}^N)^*$.

For $\mu \in C_b(\mathbb{R}^N)^*$, we define

$$I(\mu) = \sup_{V \in C_b(\mathbb{R}^N)}\{\langle V, \mu\rangle - \Lambda^*(V)\}.$$

PROPOSITION 5.4. *Let $V$ be bounded continuous. Then*

$$\Lambda^*(V) = \sup_{\mu \in C_b(\mathbb{R}^N)^*}\{\langle V, \mu\rangle - I(\mu)\}.$$

See [7], Proposition 4.1, Chapter 1 or [23], Theorem 7.15.

We shall prove that $I(\mu) = J(\mu)$ if $\mu$ is a probability measure. For $V$ satisfying (5.2), the supremum is attained at $\mu = \mu^*$, where $\mu^*$ is the invariant measure of the diffusion in (3.1) for $W = W^*$. Our main theorem is a consequence of this. We begin with some elementary observations. We follow essentially the argument in [23].

LEMMA 5.5. *If $I(\mu) < \infty$, then $\mu$ is nonnegative and $\mu(\mathbb{R}^N) = 1$.*

PROOF.    We first prove $I(\mu) = \infty$ if $\mu$ is not nonnegative. For such $\mu$, we take $V \geq 0$ such that $\mu(V) < 0$. For any $\alpha > 0$,

$$I(\mu) \geq \langle -\alpha V, \mu \rangle - \Lambda^*(-\alpha V).$$

Now $\Lambda^*(-\alpha V) \leq \Lambda^*(0)$, since $-\alpha V \leq 0$. Therefore,

$$I(\mu) \geq -\alpha \langle V, \mu \rangle - \Lambda^*(0) \to \infty$$

as $\alpha$ tends to infinity. Hence $I(\mu) = \infty$.

We now prove $I(\mu) = \infty$ if $\mu(\mathbb{R}^N) \neq 1$. For such $\mu$,

$$I(\mu) \geq \langle \alpha, \mu \rangle - \Lambda^*(\alpha) = \alpha\mu(\mathbb{R}^N) - \alpha - \Lambda^*(0) = \alpha(\mu(\mathbb{R}^N) - 1) - \Lambda^*(0).$$

Here $\alpha$ is any real number. From this, it is easy to see that $I(\mu) = \infty$.  □

LEMMA 5.6.    *Let $\mu$ be a probability measure. Then $I(\mu) \geq J(\mu)$.*

PROOF.    Let $W$ be a smooth function such that $G(W)$ is bounded above. Take

$$V_n = \min\{-G(W), n\}.$$

Then $V_n$ is a bounded continuous function.

It is easy to see that

$$\tfrac{1}{2}D_i(a^{ij}D_jW) + \tfrac{1}{2}\hat{a}^{ij}D_iWD_jW + b \cdot \nabla W + V_0 + V_n \leq 0.$$

Therefore, $\Lambda^*(V_n) \leq 0$. From the relation

$$I(\mu) \geq \langle V_n, \mu \rangle - \Lambda^*(V_n) \geq \langle V_n, \mu \rangle,$$

and $V_n \to -G(W)$, $V_n \leq V_{n+1}$ such that $V_n$ are bounded below, we can apply the monotone convergence theorem to get

$$I(\mu) \geq -\langle G(W), \mu \rangle.$$

Then

$$I(\mu) \geq \sup\{-\langle G(W), \mu \rangle; W \text{ is smooth}, G(W) \text{ is bounded above}\} = J(\mu).□$$

LEMMA 5.7.    *Let $V \in C_b(\mathbb{R}^N)$. Define*

$$J^*(V) = \sup_{\mu \in \mathcal{M}_1(\mathbb{R}^N)} \{\langle V, \mu \rangle - J(\mu)\}.$$

*Then $J^*(V) \leq \Lambda^*(V)$.*

PROOF.    It is enough to prove this for $V$ satisfying (5.2) which we assume now. We first observe the relation

$$
\begin{aligned}
J^*(V) &= \sup_\mu \left\{ \langle V, \mu \rangle + \inf_W \{ \langle G(W), \mu \rangle \} \right\} \\
&= \sup_\mu \inf_W \{ \langle G(W) + V, \mu \rangle \} \\
&\leq \inf_W \sup_\mu \{ \langle G(W) + V, \mu \rangle \} \\
&= \inf_W \sup_{x \in \mathbb{R}^N} \{ G(W)(x) + V(x) \}.
\end{aligned}
$$

Here the $\mu, W$ are taken over all those satisfying $\mu \in \mathcal{M}_1(\mathbb{R}^N)$ and $W$ smooth with $G(W)$ bounded above. Let $W^*$ be the solution of (5.1) for $\Lambda = \Lambda^*(V)$. Take $W = W^*$ in the above relation; we have $J^*(V) \leq \Lambda^*(V)$.   $\square$

LEMMA 5.8.    *Let $V$ satisfy* (5.2) *and let $\mu^*$ be the invariant measure of the diffusion in* (3.1) *for $\Lambda = \Lambda^*(V)$ and $W^*$ the solution of* (5.1). *Then*

$$
\Lambda^*(V) = \langle V, \mu^* \rangle - I(\mu^*).
$$

PROOF.    We need to prove that

$$
\langle V', \mu^* \rangle - \Lambda^*(V') \leq \langle V, \mu^* \rangle - \Lambda^*(V)
$$

for all $V' \in C_b(\mathbb{R}^N)$. This is equivalent to

$$
\Lambda^*(V') - \Lambda^*(V) \geq \langle V' - V, \mu^* \rangle, \qquad V' \in C_b(\mathbb{R}^N).
$$

Since $\Lambda^*(\cdot)$ is a convex function on $C_b(\mathbb{R}^N)$, there is a subgradient $\bar{\mu} \in C_b(\mathbb{R}^N)^*$ of this function at $V$ such that

$$
\Lambda^*(V') - \Lambda^*(V) \geq \langle V' - V, \bar{\mu} \rangle, \qquad V' \in C_b(\mathbb{R}^N).
$$

See [7], Proposition 5.2, Chapter 1. We only need to prove $\bar{\mu} = \mu^*$, since this implies the claim.

First, the nondecreasing of $\Lambda^*(\cdot)$ implies $\bar{\mu}$ is nonnegative.

Applying the above relation to $V' = V + \alpha$, $\alpha$ is constant, and using $\Lambda^*(V + \alpha) = \Lambda^*(V) + \alpha$, we can easily deduce $\bar{\mu}(\mathbb{R}^N) = 1$.

Now take $\phi_n$ smooth functions on $\mathbb{R}^N$ satisfying the following properties: $0 \leq \phi_n \leq \phi_{n+1} \leq 1$, $\phi_n$ has compact support, $\nabla \phi_n$ are bounded uniformly in $n$, $\phi_n \to 1$ uniformly on compact sets. Then Proposition 4.1 implies $\Lambda^*(V + \alpha \phi_n) \to \Lambda^*(V + \alpha) = \Lambda^*(V) + \alpha$ as $n \to \infty$. Then

$$
\limsup_{n \to \infty} \alpha \langle \phi_n, \bar{\mu} \rangle \leq \alpha
$$

for all $\alpha$. From this, we have

$$
\lim_{n \to \infty} \langle 1 - \phi_n, \bar{\mu} \rangle = 0.
$$

Then $\bar{\mu}$ must be a probability measure. This is a consequence of [6], Theorem III.5.13.

We now prove

$$\int (\tfrac{1}{2}a^{ij}(x)D_{ij}W(x) + (b(x) + \hat{a}\nabla W^*(x)) \cdot \nabla W(x))\,d\bar{\mu}(x) = 0$$

for $W$ smooth function on $\mathbb{R}^n$ with compact support. This implies $\bar{\mu}$ is an invariant measure of the diffusion in (3.1). Hence $\bar{\mu} = \mu^*$ by the uniqueness of the invariance measure. To prove this last statement, we take such $W$ and consider

$$V' = V - (\tfrac{1}{2}a^{ij}D_{ij}W + (b + \hat{a}\nabla W^*) \cdot \nabla W + \tfrac{1}{2}\hat{a}^{ij}D_i W D_j W).$$

Then by a simple calculation, we see

$$\tfrac{1}{2}a^{ij}D_{ij}(W + W^*) + b \cdot \nabla(W + W^*)$$
$$+ \tfrac{1}{2}\hat{a}^{ij}D_i(W + W^*)D_j(W + W^*) + V_0 + V' = \Lambda^*(V).$$

Therefore, $\Lambda^*(V') = \Lambda^*(V)$. Then,

$$0 \geq \langle V' - V, \bar{\mu}\rangle = -\langle \tfrac{1}{2}a^{ij}D_{ij}W + (b + \hat{a}\nabla W^*) \cdot \nabla W + \tfrac{1}{2}\hat{a}^{ij}D_i W D_j W, \bar{\mu}\rangle.$$

We replace $W$ by $\alpha W$, $\alpha > 0$, divide the relation by $\alpha$ and let $\alpha$ tend to 0. We get

$$-\langle \tfrac{1}{2}a^{ij}D_{ij}W + (b + \hat{a}\nabla W^*) \cdot \nabla W, \bar{\mu}\rangle \leq 0.$$

We replace $W$ by $-W$. Then we find

$$\langle \tfrac{1}{2}a^{ij}D_{ij}W + (b + \hat{a}\nabla W^*)\nabla W, \bar{\mu}\rangle = 0.$$

This is what we want to prove. $\quad\square$

COROLLARY 5.9.  *Let $V \in C_b(\mathbb{R}^N)$. Then $J^*(V) = \Lambda^*(V)$. Let $\mu$ be a probability measure. Then $I(\mu) = J(\mu)$.*

PROOF.  We assume $V \in C_b(\mathbb{R}^N)$ satisfying (5.2). Let $\mu^*$ be the invariant measure for the diffusion in (3.1) with $W = W^*$. Then

$$J^*(V) \geq \langle V, \mu^*\rangle - J(\mu^*) \geq \langle V, \mu^*\rangle - I(\mu^*) = \Lambda^*(V).$$

But we have already proved $J^*(V) \leq \Lambda^*(V)$. Therefore, they are equal.

We now prove $I(\mu) = J(\mu)$. We use the relation

$$J(\mu) = \sup_{V \in C_b(\mathbb{R}^N)} \{\langle V, \mu\rangle - J^*(V)\}.$$

See [23], Theorem 7.18. By definition,

$$I(\mu) = \sup_{V \in C_b(\mathbb{R}^N)} \{\langle V, \mu\rangle - \Lambda^*(V)\}.$$

Since $J^*(V) = \Lambda^*(V)$ for all $V$, we have $I(\mu) = J(\mu)$. $\quad\square$

**Acknowledgments.** The authors would like to thank the referees for helpful comments and suggestions. This research was done while the first author was a post-doctoral Research Fellow at Institute of Mathematics, Academia Sinica, Taiwan. The first author is grateful to Academia Sinica for giving him an opportunity to visit there for his research.

## REFERENCES

[1] Bensoussan, A. (1988). *Perturbation Methods in Optimal Control.* Wiley, New York. MR0949208

[2] Bensoussan, A. and Frehse, J. (1992). On Bellman equations of ergodic control in $\mathbb{R}^N$. *J. Reine Angew. Math.* **429** 125–160. MR1173120

[3] Bielecki, T. R. and Pliska, S. R. (1999). Risk sensitive dynamic asset management. *Appl. Math. Optim.* **39** 337–360. MR1675114

[4] Donsker, M. D. and Varadhan, S. R. S. (1975). Asymptotic evaluation of certain Wiener integrals for large time. In *Functional Integration and Its Applications* (A. M. Arthurs, ed.) 15–33. Oxford Univ. Press. MR0486395

[5] Donsker, M. D. and Varadhan, S. R. S. (1976). Asymptotic evaluation of certain Markov processes expectations for large time III. *Comm. Pure Appl. Math.* **29** 389–461. MR0428471

[6] Dunford, N. and Schwartz, J. T. (1988). *Linear Operators Part I: General Theory.* Wiley, New York. MR1009162

[7] Ekeland, I. and Teman, R. (1976). *Convex Analysis and Variational Problems.* North-Holland, Amsterdam. MR0463994

[8] Fleming, W. H. (1995). Optimal investment model and risk-sensitive stochastic control. In *IMA Vols. in Math. Appl.* **65** 75–88. Springer, New York.

[9] Fleming, W. H. and James, M. R. (1995). The risk-sensitive index and the $H_2$ and $H_\infty$ norms for nonlinear systems. *Math. Control Signals Systems* **8** 199–221. MR1387043

[10] Fleming, W. H. and McEneaney, W. M. (1995). Risk-sensitive control on an infinite time horizon. *SIAM J. Control Optim.* **33** 1881–1915. MR1358100

[11] Fleming, W. H. and Sheu, S.-J. (1999). Optimal long term growth rate of expected utility of wealth. *Ann. Appl. Probab.* **9** 871–903. MR1722286

[12] Fleming, W. H. and Sheu, S.-J. (2000). Risk sensitive control and an optimal investment model. *Math. Finance* **10** 197–213. MR1802598

[13] Kaise, H. and Nagai, H. (1998). Bellman–Isaacs equations of ergodic type related to risk-sensitive control and their singular limits. *Asymptotic Anal.* **16** 347–362. MR1612829

[14] Kaise, H. and Nagai, H. (1999). Ergodic type Bellman equations of risk-sensitive control with large parameters and their singular limits. *Asymptotic Anal.* **20** 279–299. MR1715337

[15] Kaise, H. and Sheu, S.-J. (2004). Risk-sensitive optimal investment: Solutions of dynamical programming equation. In *Mathematics of Finance, Contemporary Math.* **351**. Amer. Math. Soc., Providence, RI. MR2076543

[16] Kaise, H. and Sheu, S.-J. (2004). On the structure of solutions of ergodic type Bellman equation related to risk-sensitive control. Technical report, Academia Sinica.

[17] Kaise, H. and Sheu, S.-J. (2004). Evaluation of large time expectations for diffusion processes. Technical report, Academia Sinica.

[18] Ladyzhenskaya, O. A. and Ural'tseva, N. N. (1968). *Linear and Quasilinear Elliptic Equations*. Academic Press, New York. MR0244627

[19] McEneaney, W. M. and Ito, K. (1997). Infinite time-horizon risk-sensitive systems with quadratic growth. In *Proceedings of 36th IEEE Conference on Decision and Control*.

[20] Nagai, H. (1996). Bellman equation of risk-sensitive control. *SIAM J. Control Optim.* **34** 74–101. MR1372906

[21] Nagai, H. (2003). Optimal strategies for risk-sensitive portfolio optimization problems for general factor models. *SIAM J. Control Optim.* **41** 1779–1800. MR1972534

[22] Pinsky, R. G. (1995). *Positive Harmonic Functions and Diffusion*. Cambridge Univ. Press. MR1326606

[23] Stroock, D. W. (1984). *An Introduction to Theory of Large Deviations*. Springer, Berlin. MR0755154

[24] Varadhan, S. R. S. (1980). *Diffusion Problems and Partial Differential Equations*. Springer, New York.

Graduate School of Information Science          Institute of Mathematics
Nagoya University                               Academia Sinica
Furo-cho, Chikusa-ku                            Nankang, Taipei 11529
Nagoya 464-8601                                 Taiwan
Japan                                           Republic of China
e-mail: kaise@is.nagoya-u.ac.jp                 e-mail: sheusj@math.sinica.edu.tw